\documentclass[10pt]{article}

\date{\empty}

\usepackage{amsthm}
\usepackage[normalem]{ulem} 
\newtheorem{theorem}{Theorem}[section]

\newtheorem{defn}[theorem]{Definition}
\usepackage{nicematrix}

\theoremstyle{remark}
\newtheorem{remark}{Remark}[section]

\newcommand{\R}{\mathbb{R}}

\newcommand{\A}{\mathcal{A}}

\newcommand{\tpose}{^\mathsf{T}}

\DeclareMathOperator*{\argmin}{argmin}
\DeclareMathOperator*{\argmax}{argmax}

\DeclareMathOperator{\dist}{dist}

\definecolor{myblue}{rgb}{0.21, 0.34, 0.74}
\definecolor{myred}{rgb}{0.79, 0.0, 0.09}
\definecolor{mygreen}{rgb}{0, 0.32, 0}

\usepackage{amssymb}
\usepackage{amsmath}
\usepackage{graphicx}
\usepackage{cite}
\usepackage{xcolor}
\usepackage{hyperref}
\hypersetup{
    colorlinks=true,
    linkcolor=blue,
    filecolor=magenta,      
    urlcolor=cyan,
    }
\usepackage{bigints}
\usepackage{booktabs}
\usepackage{bm}
\usepackage{amsmath}
\usepackage{amsthm}
\usepackage{amssymb}
\usepackage{algorithm,algpseudocode}
\algblock{ParFor}{EndParFor}
\algnewcommand\algorithmicparfor{\textbf{parfor}}
\algnewcommand\algorithmicpardo{\textbf{do}}
\algnewcommand\algorithmicendparfor{\textbf{end\ parfor}}
\algrenewtext{ParFor}[1]{\algorithmicparfor\ #1\ \algorithmicpardo}
\algrenewtext{EndParFor}{\algorithmicendparfor}
\usepackage{amsfonts}
\usepackage{mathtools}
\usepackage{dsfont}
\theoremstyle{definition}
\newtheorem{thm}{Theorem}
\theoremstyle{remark}

\usepackage{enumitem}
\usepackage{bbm}
\usepackage{booktabs}
\usepackage{soul}
\DeclareMathOperator{\diag}{diag}
\DeclareMathOperator{\Diag}{Diag}
\usepackage{tikz}

\usepackage{cancel}

\DeclareMathOperator{\tr}{tr}
\usepackage[version=4]{mhchem}
\usepackage{graphicx}
\usepackage{fontawesome5}
\usepackage{float}

\usepackage{tcolorbox} 
\newtcolorbox{myabstract}{colback=gray!5, colframe=gray!5, sharp corners=all, boxrule=2mm}

\usepackage{fancyhdr} 
\begin{document}
\title{
{\Large \bf Modified Patankar Semi-Lagrangian Scheme for the Optimal Control of Production-Destruction systems}
}
\vspace{1em}
\author{
{\large Simone Cacace}{\small $^{1,3}$} \href{https://orcid.org/0000-0002-5864-0970}{\includegraphics[width=0.9em]{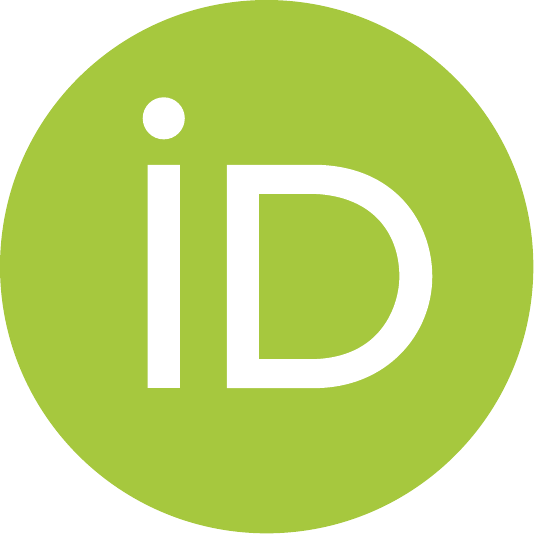}} \and
{\large Alessio Oliviero}{\small $^{1,3}$} \href{https://orcid.org/0009-0003-2603-6431}{\includegraphics[width=0.9em]{figures/orcid.pdf}} \and
{\large Mario Pezzella}{\small$^{2,3,4}$} \href{https://orcid.org/0000-0002-1869-945X}{\includegraphics[width=0.9em]{figures/orcid.pdf}}
}

\footnotetext[1]{\footnotesize Sapienza University of Rome, Department of Mathematics ``Guido Castelnuovo'', Piazzale Aldo Moro, 5 - 00185 Rome, Italy.}
\footnotetext[2]{\footnotesize C.N.R. National Research Council of Italy, Institute for Applied Mathematics ``Mauro Picone'', Via P. Castellino, 111 - 80131 Naples, Italy.} 
\footnotetext[3]{\footnotesize Member of the INdAM Research group GNCS}
\footnotetext[4]{\footnotesize Corresponding Author \href{mailto:mario.pezzella@cnr.it}{\faEnvelope[regular]}}

\date{}
\maketitle

\begin{myabstract}
\begin{center}
\textbf{Abstract}
\end{center}
\vspace{-1em}
In this manuscript, we present a comprehensive theoretical and numerical framework for the control of production-destruction differential systems. The general finite horizon optimal control problem is formulated and addressed through the dynamic programming approach. We develop a parallel in space conservative scheme for the corresponding backward-in-time Hamilton-Jacobi-Bellman equation. Furthermore, we provide a suitable reconstruction algorithm for optimal controls and trajectories. The application to two case studies, specifically enzyme catalyzed biochemical reactions and infectious diseases, highlights the advantages of the proposed methodology over classical semi-Lagrangian discretizations.\\
\phantom{ZZz} 
\end{myabstract}

\noindent\rule{\linewidth}{0.5pt} 

\section{Introduction}\label{sec:introduction}
The starting point of this work is the following Production-Destruction differential System (PDS)
\begin{equation}\label{eq:PDS_not_controlled}
\left\{
\begin{aligned}
    \bm{y}'(t)&=(P(\bm{y}(t))-D(\bm{y}(t))) \bm{e}, \\
    \bm{y}(t_0)&=\bm{y}^0, \qquad \qquad \qquad  \qquad \qquad \qquad  \qquad \qquad  t\geq t_0\geq 0,
\end{aligned}
\right.
\end{equation}
where $\bm{y}(t)=(y_1(t),\dots,y_N(t))\tpose{} \in \mathbb{R}^N,$ with $N>1,$ represents the state variables vector, $\bm{e}=(1,\dots,1)\tpose{}\in \mathbb{R}^N$ is a vector of all ones and the matrices $P(\bm{y}(t))=\{P_{ij}(\bm{y}(t))\}\in \mathbb{R}^{N\times N}$ and $D(\bm{y}(t))=\{D_{ij}(\bm{y}(t))\}\in \mathbb{R}^{N\times N}$ denote the non-linear production and destruction rates, respectively. Numerous real-world phenomena involve the interplay between processes of production and decay or consumption and can be therefore modeled by PDS of the form \eqref{eq:PDS_not_controlled} (see, for instance, \cite{Higham,CdS,Epidemic_PDS_1,Epidemic_PDS_2,Bacteria_PDS}.)

Our investigation will be based on the following assumptions for the system \eqref{eq:PDS_not_controlled}
\begin{enumerate}[label=(A$_{\arabic*}$)]
    \item\label{ass:A1} $P_{ij}(\bm{x})$ and $D_{ij}(\bm{x})$ are non-negative, continuous functions on $\mathbb{R}_{0^+}^N;$  
    \item\label{ass:A2} $P(\bm{x})=D\tpose(\bm{x}),$ for all $\bm{x}\in \mathbb{R}^{N}_{0^+};$
    \item\label{ass:A3} $P_{ii}(\bm{x})=D_{ii}(\bm{x})=0,$ for all $\bm{x}\in \mathbb{R}^{N}_{0^+};$
    \item\label{ass:A4} $P_{ij}(x_1,\dots,x_{j-1},0,x_{j+1},\dots,x_N)=0,$ for all $(x_1,\dots,x_{j-1},x_{j+1},\dots,x_N)\in \mathbb{R}^{N-1}_{0^+};$
\end{enumerate} 
for all $i,j=1,\dots,N$, with
\begin{equation*}
     \mathbb{R}^{N}_{0^+}=\left\{\bm{x}=(x_1,\dots,x_N)\tpose\in\mathbb{R}^N \ : \ x_i\geq 0, \;\; \mbox{for} \ i=1,\dots,N\right\}.
\end{equation*}
The properties above ensure the existence of a unique non-negative solution to the initial value problem \eqref{eq:PDS_not_controlled} (see \cite[Theorem 3.3]{Formaggia_Scotti} and \cite[Theorem 1.2]{TORLO2022} for further discussion on this topic). Furthermore, owing to \ref{ass:A2} and \ref{ass:A3}, the PDS \eqref{eq:PDS_not_controlled} is positive and fully conservative (see \cite[Definitions 1.1-1.2]{Kopecz2018}), i.e.
\begin{equation}\label{eq:Omega0}
    \Omega_0=\{\bm{x}\in\mathbb{R}^N \ : \ \bm{0}\leq \bm{x} \leq \bm{e}\tpose \bm{y}^0  \},
\end{equation}
is a positively invariant set for the system \eqref{eq:PDS_not_controlled}, whose solution satisfies
\begin{equation}\label{eq:conservation_law_PDS}
    \bm{y}^0>0 \implies \bm{y}(t)>0, \quad \forall t \geq t_0, \qquad \mbox{and} \qquad
    \bm{e}\tpose\bm{y}(t)=\bm{e}\tpose\bm{y}^0, \quad \forall t \geq t_0.
\end{equation}
In \eqref{eq:conservation_law_PDS} and throughout the paper, inequalities involving vectors are considered component-wise. 

The numerical simulation of positive and fully conservative PDS presents challenges in retaining the properties outlined in \eqref{eq:conservation_law_PDS}, which are guaranteed by standard methods only for sufficiently small time steps. Therefore, the necessity arises for unconditionally positive and conservative schemes that, assuming positive initial values, yield positive numerical solutions and preserve the linear invariant of system \eqref{eq:PDS_not_controlled}, independently of the discretization step-length. Numerous contributions in the literature have addressed this topic (see, for instance, \cite{Dimitrov1,Dimitrov2,GECO1,GECO2,Formaggia_Scotti}). Of particular relevance are the Modified Patankar methods, linearly implicit time integrators devised by suitably modifying explicit schemes. The original approach, introduced by Patankar in \cite{patankar1980numerical}, has been firstly extended to Runge-Kutta discretizations \cite{MPEuler,Kopecz2018,Kopecz2018second,Kopecz2019,Huang2019,Huang2019II,Izgin3,Izgin1,IzginSSPMPRK,Izgin2} and more recently applied to deferred correction \cite{TORLO2022,Torlo2020} and linear multistep \cite{MPLM,Zhu2024} methods . 

Optimal control theory provides a valuable framework to guide interventions across the variety of natural and industrial processes modeled by PDS. The approaches to address optimal control problems can be classified into two categories. The direct methods, such as the forward-backward sweep or the direct-adjoint looping \cite{MMH12}, rely on Pontryagin's first order optimality conditions \cite{BGP56}. Alternatively, indirect methods involve a reformulation of the problem in terms of partial differential equations \cite{BardiCapuzzoDocetta2008,FalconeFerretti}. More recently, hybrid methods that combine the aforementioned procedures \cite{CM10,AKK21,Cacace_Oliviero} and machine learning techniques \cite{SSGK23} have become viable thanks to the advancements in GPUs utilization.

To the best of our knowledge, there are no contributions in the literature that specifically pertain to the optimal control of PDS. However, the versatility of these models and the need to preserve their inherent properties motivate us in developing a theoretical and numerical framework for this purpose, which we present in this work.
We adopt the dynamic programming \cite{Bellman1957} indirect approach that guarantees convergence to the optimal solution and provides optimal controls in \textit{feedback} form. First, from the optimal control problem a first-order Hamilton-Jacobi equation for the value function is formulated. Then, the optimal trajectory and control are synthesized by solving a finite-dimensional optimization problem. The most challenging aspect of this routine lies on the approximation of the PDE viscosity solution \cite{CrandallLions1983,CrandallEvansLions1984}. As a matter of fact, commonly employed numerical methods such as finite differences may reveal inadequate for the task, as they rely on regularity assumptions for the solution. In this context, semi-Lagrangian schemes represent a valid alternative, providing greater stability and less numerical diffusion by tracking the evolution of particles along characteristic paths. Here, we define a robust and reliable semi-Lagrangian method built upon the approximation of the PDS characteristics via modified Patankar integrators, taking advantage of their conservativity and positivity properties. The proposed scheme is, by design, parallelizable in space, which allows for an efficient implementation on GPUs.

The manuscript is structured as follows. Section \ref{sec:CPS} introduces the theoretical framework for the controlled production-destruction system. There, the optimal control problem and the dynamic programming approach are detailed as well. In Section \ref{sec:MPSL}, we present the parallel-in-space modified Patankar semi-Lagrangian scheme and a positive and conservative algorithm for the reconstruction of the optimal control and trajectory. Section \ref{sec:numerical} addresses two specific case studies, supported by numerical simulations. Finally, Section \ref{sec:Conclusions} offers concluding remarks and outlines potential directions for future research.

\section{Controlled production-destruction systems}\label{sec:CPS} 

Let $M\leq N$ be a positive integer, $A$ be a compact subset of $\mathbb{R}^M$ and
\begin{equation}\label{eq:control_function}
    \bm{\alpha}: \; t\in [0,+\infty) \longrightarrow \bm{\alpha}(t)\in A, \; \qquad \; \mathcal{P},\mathcal{D}: \bm{a}\in A \longrightarrow \mathcal{P}(\bm{a}),\mathcal{D}(\bm{a})\in \mathbb{R}^{N\times N},
\end{equation}
be Lebesgue-measurable functions. From now on, we assume that $\bm{y}^0 \geq 0$ and that \ref{ass:A1}-\ref{ass:A4} hold true. Referring to the notations in Section \ref{sec:introduction}, we then define a Controlled Production-Destruction System (CPDS) as follows
\begin{equation}\label{eq:CPDS}
\left\{
\begin{aligned}
    \bm{y}'(t)&=\big(P(\bm{y}(t))\odot \mathcal{P}(\bm{\alpha}(t))-D(\bm{y}(t))\odot \mathcal{D}(\bm{\alpha}(t)) \big) \bm{e}, \\
    \bm{y}(t_0)&=\bm{y}^0, \qquad \qquad \qquad  \qquad \qquad \qquad  \qquad \qquad \qquad \qquad  t\geq t_0\geq 0,
\end{aligned}
\right.
\end{equation}
where the symbol $\odot$ denotes the (component-wise) Hadamard product. In this context, $\bm{\alpha}(t)$ represents the control function at time $t\geq t_0.$ The matrices $\mathcal{P}(\bm{\alpha}(t))=\{\mathcal{P}_{ij}(\bm{\alpha}(t))\}\in\mathbb{R}^{N\times N}$ and $\mathcal{D}(\bm{\alpha}(t))=\{\mathcal{D}_{ij}(\bm{\alpha}(t))\}\in\mathbb{R}^{N\times N}$ govern the control policies  associated with the production and destruction processes, respectively. The following definitions arise from the aim of controlling a PDS while avoiding any alterations to its inherent nature and properties.

\begin{defn}\label{defn:Positive_CPDS}
    The controlled system \eqref{eq:CPDS} is referred to as positive if, independently of the choice of the control function in the space of admissible controls
    \begin{equation}\label{eq:set_admissible_controls}
        \mathcal{A}=\{\bm{\alpha}: [0,+\infty) \longrightarrow A \;\; | \ \bm{\alpha} \; \mbox{measurable} \},
    \end{equation}
    the condition $\bm{y}^0>0$ implies $\bm{y}(t)>0$ for all $t\geq t_0.$
\end{defn}

\begin{defn}\label{defn:Conservativity_CPDS}
    The controlled system \eqref{eq:CPDS} is called conservative if for all $\bm{\alpha}\in\mathcal{A},$ the following conservation law holds
    \begin{equation*}
        \bm{e}\tpose(\bm{y}(t)-\bm{y}^0)=0, \qquad \mbox{for all } t\geq t_0.
    \end{equation*}
\end{defn}

The existence of a unique solution to the CPDS \eqref{eq:CPDS} is proved, under mild regularity assumptions on the known functions, with the following result.

\begin{thm}[\textbf{Existence of the solution}]\label{thm:Existence of a solution}
    Assume that $P$ and $D$ are Lipschitz continuous functions on $\mathbb{R}^N$ and that $\mathcal{P},\mathcal{D}\in C^0(A).$  Then, for each $\bm{\alpha}\in\mathcal{A},$ the controlled production-destruction system \eqref{eq:CPDS} admits a unique solution. 
\end{thm}
\begin{proof}
    The CPDS \eqref{eq:CPDS} can be rewritten as $\bm{y}'(t)= \bm{F}(\bm{y}(t),\bm{\alpha}(t))$, with $\bm{y}(0)=\bm{y}^0,$  $\bm{\alpha}\in \mathcal{A}$ and
    \begin{equation}\label{eq:F}
        \bm{F} : (\bm{\mu},\bm{a})\in \mathbb{R}^N\times A \longrightarrow \big(P(\bm{\mu})\odot \mathcal{P}(\bm{a})-D(\bm{\mu})\odot \mathcal{D}(\bm{a}) \big) \bm{e}  \in \mathbb{R}^N.
    \end{equation}
    Let $L_P,$ $L_D$ be the Lipschitz constants of $P$ and $D$, respectively. Define 
    \begin{equation*}
        \bar{p}=\!\!\max_{\substack{\bm{a}\in A \\ 1\leq i,j\leq N}}|\mathcal{P}_{ij}(\bm{a})| \qquad \qquad \text{and} \qquad \qquad  \bar{d}=\!\!\max_{\substack{\bm{a}\in A \\ 1\leq i,j\leq N}}|\mathcal{D}_{ij}(\bm{a})|.
    \end{equation*}
    The hypotheses on the known functions imply that $\bm{F}\in C^0(\mathbb{R}^N\times A)$ and that, for each $\bm{\mu},\bm{\nu} \in \mathbb{R}^N$ and $\bm{a}\in A,$ 
    \begin{equation*}
        \begin{split}
            \|\bm{F}(\bm{\mu},\bm{a})-\bm{F}(\bm{\nu},\bm{a}) \| \leq \; & \hat{e}\|\big (P(\bm{\mu})-P(\bm{\nu}))\odot \mathcal{P}(\bm{a})  \|+\hat{e}\|\big (D(\bm{\mu})-D(\bm{\nu}))\odot \mathcal{D}(\bm{a}) \| \\
            \leq \; &\hat{e}k_N\left(\bar{p}L_P  +\bar{d}L_D\right)\|\bm{\mu}-\bm{\nu}\|,
        \end{split}
    \end{equation*}
    where $\hat{e}=\|\bm{e}\|$ and $k_N>0$ is a constant depending on $N.$ As a result, the function $\bm{F}$ in \eqref{eq:F} is Lipschitz continuous as well and all the hypotheses of the Carathéodory theorem (see \cite[Theorem 8.1]{FalconeFerretti}) are fulfilled. The existence of a unique solution $\bm{y}(t)$ to \eqref{eq:CPDS} is then assured.
\end{proof}

\begin{remark}
    The solution given by Theorem \ref{thm:Existence of a solution} must be regarded as a weak, almost everywhere differentiable, solution of the system \eqref{eq:CPDS}.
\end{remark}

To establish physically meaningful controls policies, we investigate sufficient conditions for the positivity and conservativity of a CPDS.

\begin{thm}[\textbf{Positivity of CPDS}]\label{thm:Positive_CPDS}
    Assume that $P$ and $D$ are Lipschitz continuous functions on $\mathbb{R}^N$ and that $\mathcal{P},\mathcal{D}\in C^0(A).$ Then the CPDS  \eqref{eq:CPDS} is positive, according to the Definition \ref{defn:Positive_CPDS}.
\end{thm}
\begin{proof}
    We emphasize that, due to \ref{ass:A4}, the condition $\bm{y}^0=\bm{0}$ implies  $\bm{y}(t)=\bm{0}$ for all $t\geq t_0.$ We prove the positivity of the CPDS by contradiction. Let $\bm{y}^0>\bm{0}$ and suppose that there exist $i^*\in \{1,\dots,N\}$ and $t^*>t_0,$ such that $y_{i^*}(t^*)\leq 0.$ It then follows that the non-empty set 
    \begin{equation*}
        \mathfrak{T}=\{\tau \geq t_0 \ : \ \bm{y}(t)>\bm{0} \quad \mbox{for all} \ t\in [t_0,\tau] \},
    \end{equation*}
    is bounded from above and that $\bar{t}=\sup \mathfrak{T}\in [t_0,t^*).$ We define, for $\varepsilon >0,$ the perturbed controlled system 
    \begin{equation}\label{eq:shifted_CPDS}
    \left\{
    \begin{aligned}
        \bm{y}^{\varepsilon \prime}(t)&=\bm{F}(\bm{y}^{\varepsilon}(t),\bm{\alpha}(t)) + \varepsilon \bm{e}, \\
        \bm{y}^{\varepsilon}(t_0)&=\bm{y}^0+ \varepsilon \bm{e},
    \end{aligned}
    \right.
    \end{equation}
    with $\bm{F}(\bm{\mu},\bm{a})=(F_1(\bm{\mu},\bm{a}),\dots,F_N(\bm{\mu},\bm{a}))\tpose$ given in \eqref{eq:F} and denote by $\bm{y}^{\varepsilon}(t)$ its unique solution for $t\geq t_0,$ which exists for a straightforward extension of Theorem \ref{thm:Existence of a solution}. From the continuity of $\bm{F}$ and the uniqueness of the solution to \eqref{eq:CPDS}, $\lim_{\varepsilon \to 0^+}\bm{y}^{\varepsilon}(t)=\bm{y}(t)$ (see, for instance, \cite[p. 18]{coppel1965stability}). It is then possible to select $\varepsilon^*>0$ and $\tilde{t}\in (t_0,\bar{t})$ such that ${y}^{\varepsilon^*}_{i^*}(\tilde{t})=0$ and ${y}^{\varepsilon^*}_{i^*}(t)>0,$ for all $t\in [t_0,\tilde{t}).$ Therefore ${y}^{\varepsilon^* \prime}_{i^*}(\tilde{t})\leq 0.$ On the other hand, from \eqref{eq:shifted_CPDS}, \ref{ass:A2} and \ref{ass:A4}, 
    \begin{equation*}
        {y}^{\varepsilon^* \prime}_{i^*}(\tilde{t})= F_{i^*}(({y}^{\varepsilon^*}_1(t),\dots, y^{\varepsilon^*}_{i^*-1}(t),0,y^{\varepsilon^*}_{i^*+1}(t),\dots),\bm{\alpha}(t)) + \varepsilon^* = \varepsilon^* >0,
    \end{equation*}
    which yields the contradiction.
\end{proof}

Given a matrix $M=\{M_{ij}\}\in \mathbb{R}^{N \times N},$ we denote by $\tr(M)\in \mathbb{R}$ the trace of the matrix, defined as follows
\begin{equation*}
    \tr(M)=\bm{e}M\bm{e}\tpose=\textstyle\sum_{i=1}^NM_{ii},
\end{equation*}
and by $\diag(M)\in \mathbb{R}^{N \times N}$ the diagonal matrix whose elements are given by
\begin{equation}\label{eq:diag}
    \diag(M)_{ij}=\delta_{ij}M_{ij}, \qquad \mbox{for  } \; 1\leq i,j \leq N,
\end{equation}
where $\delta_{ij}$ represents the Kronecker delta. 

\begin{thm}[\textbf{Conservativity of CPDS}]\label{thm:Conservative_CPDS}
        Assume that the hypotheses of Theorem \ref{thm:Positive_CPDS} are fulfilled. Then, if the known functions satisfy
        \begin{equation}\label{eq:cond_cons_CPDS}
            \tr(P(\bm{x})(\mathcal{P}\tpose(\bm{a})-\mathcal{D}(\bm{a})))=0, \qquad \mbox{for all } \bm{x}\in \mathbb{R}^N \mbox{and } \bm{a}\in A,
        \end{equation}
        the CPDS \eqref{eq:CPDS} is conservative, according to the Definition \ref{defn:Conservativity_CPDS}. Moreover, $\Omega_0$ defined in \eqref{eq:Omega0} is a positively invariant set for the system.
\end{thm}
\begin{proof}
    Owing to \ref{ass:A2} and the properties of the Hadamard product of matrices (see, for instance, \cite[Lemma 7.5.2]{Horn_2012}), the equalities 
    \begin{equation*}
    \begin{split}
        \bm{e}\tpose(P(\bm{y}(t))\odot \mathcal{P}(\bm{\alpha}(t)))\bm{e}&=\bm{e}\tpose\diag(P(\bm{y}(t))\mathcal{P}\tpose(\bm{\alpha}(t)))=\tr(P(\bm{y}(t))\mathcal{P}\tpose(\bm{\alpha}(t))), \\
        \bm{e}\tpose(D(\bm{y}(t))\odot \mathcal{D}(\bm{\alpha}(t)))\bm{e}&=\bm{e}\tpose\diag(P\tpose(\bm{y}(t))\mathcal{D}\tpose(\bm{\alpha}(t)))=\tr(P(\bm{y}(t))\mathcal{D}(\bm{\alpha}(t))),
    \end{split}
    \end{equation*}
    hold true for each $t\geq t_0$ and $\alpha \in \mathcal{A}.$ It then follows from \eqref{eq:CPDS} and the hypotheses that
    \begin{equation*}
        \begin{split}
            (\bm{e}\tpose\bm{y}(t))'&=\bm{e}\tpose\big(P(\bm{y}(t))\odot \mathcal{P}(\bm{\alpha}(t))-D(\bm{y}(t))\odot \mathcal{D}(\bm{\alpha}(t)) \big) \bm{e} \\
            &= \tr(P(\bm{y}(t))(\mathcal{P}\tpose(\bm{\alpha}(t))-\mathcal{D}(\bm{\alpha}(t))))=0,
        \end{split}
    \end{equation*}
    which implies $\bm{e}\tpose\bm{y}(t)=\bm{e}\tpose\bm{y}^0$ for each $t \geq t_0.$ Thus, if $\bm{y}^0=\bm{0},$ then $\bm{y}(t)\equiv\bm{0}$ belongs to $\Omega_0=\{\bm{y}^0\}$ for all $t\geq t_0.$ On the other hand, if $\bm{y}^0>\bm{0},$ Theorem \ref{thm:Positive_CPDS} assures that $0<y_i(t)< \sum_{j=1}^Ny_j(t)=\bm{e}\tpose\bm{y}^0,$ for each $t\geq t_0$ and $i=1,\dots,N,$ which yields the result.
\end{proof}
Theorem \ref{thm:Conservative_CPDS} addresses CPDS that preserve $\bm{e}\tpose\bm{y}^0$ as a linear invariant. As a direct consequence of it, the condition $\mathcal{P}(\bm{a})=\mathcal{D}\tpose(\bm{a})$ for all $\bm{a}\in A,$ which mirrors the assumption \ref{ass:A2}, is sufficient, though in general not strictly necessary, for achieving a conservative CPDS.

\subsection{The optimal control problem}

Let $\ell: \R^N \times A \times [0,+\infty) \to \R$ and $\varphi: \R^N \to \R$ be bounded, uniformly continuous functions. Consider, for $t_f \in (t_0,+\infty)$ and $\bm{\alpha}\in \mathcal{A},$ the cost functional 
\begin{equation}\label{eq:cost_functional}
J_{\bm{y}^0,t_0}(\bm{\alpha}) = \int_{t_0}^{t_f} \ell(\bm{y}_{\bm{y}^0,t_0} (t; \bm{\alpha}),\bm{\alpha}(t),t)\ dt\ +\ \varphi(\bm{y}_{\bm{y}^0,t_0}(t_f; \bm{\alpha})),
\end{equation}
where $\bm{y}_{\bm{y}^0,t_0} (\, \cdot \, ; \bm{\alpha})$ denotes the solution to the CPDS \eqref{eq:CPDS} with initial data $\bm{y}^0,$ at the initial time $t_0,$ and implementing the control $\bm{\alpha}.$ A \textit{finite horizon} optimal control problem for the CPDS \eqref{eq:CPDS} consists in finding a control function $\bm{\alpha}^*$ in the space of admissible controls $\A,$ such that
\begin{equation} \label{eq:optimal_control}
    \bm{\alpha}^* = \argmin_{\bm{\alpha} \in \A} J_{\bm{y}^0,t_0}(\bm{\alpha}).
\end{equation}
The expression \textit{finite horizon} is here employed to stress the objective of controlling the system within a predetermined final time $t_f<+\infty$ (we refer to \cite[Section 8.1]{FalconeFerretti} for a comprehensive overview of other types of optimal control problems). The functions $\ell$ and $\varphi$ jointly contribute to the cost functional in \eqref{eq:cost_functional}. The former, called \textit{running cost} or Lagrangian, quantifies the cost associated with the dynamics throughout the entire time interval. The latter, representing the \textit{final cost}, depends solely on the state of the system at time $t_f$.

Proving the existence of a global minimum $\bm{\alpha}^*$ for the cost functional in \eqref{eq:cost_functional} is generally a challenging task. In fact, many theoretical results in the scientific literature rely on additional assumptions, primarily related to convexity, concerning both the cost functions and the dynamics. We refer to \cite[Section III.3.4]{BardiCapuzzoDocetta2008} and to \cite{Cesari1965existence, roxin1962existence, dontchev2020approximating} for further insights on the subject. In what follows we address problem \eqref{eq:optimal_control} by the dynamic programming approach which, under the hypotheses we made above, ensures the existence of at least one optimal control in feedback form without the need of assuming convexity.

\subsection{The dynamic programming approach}\label{subsec:DPP}

In this section, we present the dynamic programming approach (see, for instance, \cite{Bellman1957, BardiCapuzzoDocetta2008}, \cite[Chapter 8]{FalconeFerretti}) to finite horizon optimal control problems in the setting of CPDS. From now on, for the sake of simplicity, we denote by $(\bm{x}, t)$ the initial data and time $(\bm{y}^0, t_0)$, with a slight abuse of notation.
 
As a first step, we define the \textit{value function} $v : \R^N \times [0, t_f] \to \R$ as follows
\begin{equation} \label{eqn:value}
v(\bm{x},t) = \inf_{\bm{\alpha} \in \A} J_{\bm{x},t} (\bm{\alpha}),
\end{equation}
which represents the optimal cost of the trajectories of the CPDS \eqref{eq:CPDS} with initial conditions $\bm{y}(t)=\bm{x}$. It is well-established that the value function in \eqref{eqn:value} lacks differentiability even for linear systems and smooth cost functions (cf. \cite[Proposition 3.1]{BardiCapuzzoDocetta2008} for a regularity result). Furthermore, we highlight the dependence of $v(\bm{x},t)$ on the initial condition and the initial time in \eqref{eq:CPDS}. In particular, the spatial variable $\bm{x}$ belongs to the state-space of the dynamical system, which may not coincide with the physical space.

A fundamental result for dynamic programming is Bellman's Dynamic Programming Principle (DPP) \cite{Bellman1957, FalconeFerretti}, which in our case can be formulated as follows.
\begin{thm}[\textbf{DPP for CPDS}] \label{thm:DPP}
Under the assumptions of Theorem \ref{thm:Existence of a solution}, the value function $v$ in \eqref{eqn:value} satisfies
\begin{equation}\label{eq: DPP}
    v(\bm{x},t) = \inf_{\bm{\alpha} \in \A} \left\{ \int_t^\tau \ell(\bm{y}_{\bm{x},t} (s; \bm{\alpha}),\bm{\alpha}(s),s)\ ds + v(\bm{y}_{\bm{x},t} (\tau; \bm{\alpha}),\tau) \right\},
\end{equation}
for all $\bm{x} \in \R^N$ and $0 \leq t< \tau \leq t_f$.
\end{thm}
\begin{proof} 
Immediately comes from \cite[Proposition 3.2]{BardiCapuzzoDocetta2008} by setting a constantly null interest rate $\lambda.$
\end{proof}
Following the arguments in \cite{Bellman1957, BardiCapuzzoDocetta2008} and taking advantage of the DPP \eqref{eq: DPP}, we derive the following partial differential equation which characterizes $v$ for all $\bm{x} \in \R^N$ and $t \in [0, t_f]$
\begin{equation} \label{eq:HJB} 
\left\{
\begin{aligned}
&- \dfrac{\partial v}{\partial t}(\bm{x},t) +\max\limits_{\bm{a} \in A} \left\{ - \bm{e}\tpose \big(P(\bm{x})\odot \mathcal{P}(\bm{a}) - D(\bm{x})\odot \mathcal{D}(\bm{a}) \big)\tpose \,  \nabla v (\bm{x},t) -\ell(\bm{x},\bm{a},t) \right\}=0, \\
&v(\bm{x},t_f)=\varphi(\bm{x}).
\end{aligned}
\right.
\end{equation}

The differential problem \eqref{eq:HJB} corresponds to a backward-in-time Hamilton-Jacobi-Bellman (HJB) equation. Therefore, it may not admit a classical solution and concurrently exhibit infinitely many almost everywhere differentiable ones. Nevertheless, by adapting \cite[Section 10.3, Theorem 2]{Evans1998} to the case of CPDS, it can be proved that the value function in \eqref{eqn:value} is the unique viscosity solution to \eqref{eq:HJB}. For further details on viscosity solutions of Hamilton-Jacobi equations we refer to \cite{CrandallLions1983, CrandallEvansLions1984}.

Once $v(\bm{x},t)$ is derived from \eqref{eq:HJB}, an optimal control in feedback form satisfying \eqref{eq:optimal_control} can be obtained, at least where $v$ is smooth, by solving the optimization problem
\begin{equation}\label{eq:feedback}
\bm{\alpha}^* (\bm{x},t) = \argmax_{\bm{a} \in A} \left\{ -\bm{e}\tpose \big(P(\bm{x})\odot \mathcal{P}(\bm{a}) - D(\bm{x})\odot \mathcal{D}(\bm{a}) \big)\tpose \,  \nabla v (\bm{x},t) -\ell(\bm{x},\bm{a},t) \right\}.
\end{equation}
The interpretation of \eqref{eq:feedback} at the points where $\nabla v$ does not exist in a classical sense is rather delicate and goes beyond the scope of this work (see, for instance, \cite[Section III.2]{BardiCapuzzoDocetta2008}). Finally, the corresponding open-loop optimal control may be retrieved as
\begin{equation}\label{eq:opt_contr_no_feed}
   \tilde{\bm{\alpha}}^* (s) := \bm{\alpha}^* (\bm{y}_{\bm{y}^0,t_0} (s; \bm{\alpha}^*),s), \qquad s \geq t_0. 
\end{equation}

\begin{remark}    
Solving \eqref{eq:HJB} is the most demanding task in our approach, since \eqref{eq:feedback} is just an optimization problem over a compact subset of $\R^M$.
\end{remark}

\begin{remark}
    For practical applications, feedback controls are more convenient than open-loop ones because, as functions of the state, they can automatically adapt to small perturbations in the dynamics. 
    Another advantage of the dynamic programming approach is that if the initial data change, only \eqref{eq:feedback} needs to be solved again, without the need to compute a new value function, since $v$ is defined across the entire state-space.
\end{remark}

\section{A parallel-in-space Modified Patankar Semi-Lagrangian scheme for the HJB equation}
\label{sec:MPSL}

The HJB equation \eqref{eq:HJB} provides necessary and sufficient conditions for the optimality of a control with respect to \eqref{eq:cost_functional} (cf. \cite[Section 3.4]{BardiCapuzzoDocetta2008}). However, its inherent non-linearity poses a significant challenge in deriving its viscosity solution and therefore the necessity of employing numerical methods arises. Here, we adopt a  Semi-Lagrangian (SL) approach to device a parallel integrator which naturally tracks the characteristic curves of the HJB equation. This feature, combined with a Patankar-type unconditionally positive and conservative time-marching scheme, yields more accurate approximations of the right-hand side of \eqref{eq:HJB} and represents the primary advantage of the SL discretizations over other standard numerical methods, such as those based on finite differences. Moreover, the absence of CFL-like conditions allows for the use of larger time steps to mitigate numerical diffusion effects \cite{FalconeFerretti2016}. 

Consider $\Delta t > 0$ and a uniform partition $\{t^n\}_{0\leq n \leq \bar{n}}$ of the interval $[t_0,t_f]$, where the discrete times are denoted by $t^n=t_0+n \Delta t,$ for $n=0,\dots, \bar{n}$, and $t^{\bar{n}}=t_f$. Since from Theorem \ref{thm:Conservative_CPDS} the trajectories of the CPDS \eqref{eq:CPDS} are confined within the set $\Omega_0$ defined in \eqref{eq:Omega0}, we designate it as the computational space domain. Let $\{\bm{x}_i\}_{0\leq i \leq \bar{I}}$ represent a uniform mesh therein and denote by $V_i^{n}=V^n(\bm{x}_i)\approx V(\bm{x}_i,t^n)$ the discrete value function at time $t^n$ on the grid node $\bm{x}_i$. Following the arguments in \cite[Section 8.4.2]{FalconeFerretti}, we define, for $i=0,\dots, \bar{I}$ and $ n=\bar{n},\ldots,1$, the backward-in-time Modified Patankar Semi-Lagrangian (MPSL) scheme as follows 
\begin{equation}\label{eqn:MPSLscheme}
\left\{
    \begin{aligned}
        & V_i^{\bar{n}} = \varphi(\bm{x}_i),\\
        & M(\bm{x}_i,\bm{a})=\left(\!P(\bm{x}_i)\odot \mathcal{P}(\bm{a}) -\Diag\!\bigg(\!\left(D(\bm{x}_i)\odot \mathcal{D}(\bm{a})\right)\bm{e}\right)\!\bigg) \Diag\!\left(\!\left(\dfrac{1}{x_1},\dots,\dfrac{1}{x_N}\right)^{\!\mathsf{T}}\right), \\
        & V_i^{n-1} \! =  \min_{\bm{a} \in A} \left\{\mathcal{I} \left[ V^n\right]((I-\Delta t \ M(\bm{x}_i,\bm{a}))^{-1} \cdot \bm{x}_i)  + \Delta t \,\ell(\bm{x}_i,\bm{a},t^n)\right\},
    \end{aligned}
\right.
\end{equation}
where $\mathcal{I}[\, \cdot\, ]$ is a monotone (spatial) interpolation operator. Here, given a vector $\bm{w}\in \mathbb{R}^N,$ we denote by $\Diag\bm{w}=\diag(\bm{e}\bm{w}\tpose)\in \mathbb{R}^{N\times N},$ where the $\diag$ operator is defined in \eqref{eq:diag}. The peculiarity of the MPSL scheme \eqref{eqn:MPSLscheme} for the HJB equation, compared to classical SL approaches, lies in the use of the Modified Patankar-Euler (MPE) discretization \cite{MPEuler}  to approximate the foot of the characteristic curve. The conservative nature of this method excludes the necessity of projections of state space points onto the computational domain $\Omega_0$ (we refer to \cite{Falcone_Ferretti_High} for theoretical results on the use of general one-step time integrators). In addition, \eqref{eqn:MPSLscheme} retains the advantageous attributes of traditional SL schemes. Specifically, it avoids the explicit computation of $\nabla v$ on the whole domain by treating the discretization of the scalar product $\bm{e}\tpose \big(P(\bm{x})\odot \mathcal{P}(\bm{a}) - D(\bm{x})\odot \mathcal{D}(\bm{a}) \big)\tpose \,  \nabla v (\bm{x},t)$ in \eqref{eq:HJB} as a directional derivative.

The following result, which is adapted from \cite[Theorem 3.1]{FalconeGiorgi1999}, addresses the convergence of the numerical method \eqref{eqn:MPSLscheme}.
\begin{thm}\label{thm:SL_Convergence}
    Assume that the functions $P$ and $D$ in \eqref{eq:HJB} are Lipschitz continuous on $\mathbb{R}^N$ and that $\mathcal{P}$ and $\mathcal{D}$ are continuous on $A.$ Let $v(\bm{x},t)$ represent the exact viscosity solution of the HJB equation \eqref{eq:HJB} for $(\bm{x},t)\in \Omega_0 \times [t_0,t_f],$ with $0\leq t_0<t_f$ and $\Omega_0$ defined in \eqref{eq:Omega0}. Consider ${V^n_i}$ as the approximated solution computed by \eqref{eqn:MPSLscheme} with time step-length $\Delta t=(t_f-t_0)/\bar{n}$ and spatial mesh size $\Delta x=\|\bm{x}_{\bar{I}}-\bm{x}_{0}\|/\bar{I}.$ Then 
    \begin{equation}\label{eq:error_estimate}
        \max_{\substack{i=0,\dots, \bar{I} \\ n=0,\dots,\bar{n}}} |v(\bm{x}_i,t^n)-V_i^n|\leq C \left(\Delta t^{1/2}+\dfrac{\Delta x}{\Delta t^{1/2}}\right),
    \end{equation}
    where the positive constants $C$ is independent of $\Delta x$ and $\Delta t.$
\end{thm}

\begin{remark}
    Under the broader assumption of $\gamma$-H\"{o}lder continuity for the functions $P$ and $D$ in \eqref{eq:CPDS}, convergence of order $\gamma/2$ can be established for the MPSL scheme \eqref{eqn:MPSLscheme}.
\end{remark}

The MPSL scheme \eqref{eqn:MPSLscheme} requires the solution of an optimization problem at each time step, for which specialized algorithms can be employed. Here, to achieve a fully discrete method, we approximate the minimum through direct comparison by discretizing the compact set $A\subset \mathbb{R}^M$ with an adequate number of discrete controls. Moreover, in order to further mitigate the computational demands of \eqref{eqn:MPSLscheme}, we exploit its structure to implement a spatial parallelization procedure based on the \texttt{CUDA} programming model \cite{CudaGuide}. Specifically, although individual $V_j^{n-1}$ values may depend on the entire set $\{V_i^{n}\}_{0\leq i \leq \bar{I}}$ due to the interpolation operator $\mathcal{I}$, the computations for the different $V_i^{n-1},$ for $i=0,\dots, \bar{I},$ can occur simultaneously.

\subsection{Reconstruction of the optimal control and trajectory}\label{subsec:Opti_Traj}\label{sec:feedback_reconstruction}

The approximation of the value function provided by the MPSL method \eqref{eqn:MPSLscheme} is used to concurrently synthesize the optimal control $\tilde{\bm{\alpha}}^*$ and the optimal trajectories $\bm{y}^*$ by numerically solving \eqref{eq:feedback}-\eqref{eq:opt_contr_no_feed}. Specifically, in line with the MPSL approach, we employ the unconditionally positive and conservative MPE scheme for time integration and a direct comparison procedure for the maximization in \eqref{eq:feedback}.

Our optimal control-trajectory reconstruction algorithm then reads
\begin{equation}\label{eqn:feedback_rec_MPE}
\left\{
    \begin{aligned}
        & M(\bm{y}^{*,n}\!,\bm{a})=\Big(\!P(\bm{y}^{*,n})\odot \mathcal{P}(\bm{a}) -\Diag\!\big(\!\left(D(\bm{y}^{*,n})\odot \mathcal{D}(\bm{a})\right)\bm{e}\big)\!\Big) \Diag\!\left(\!\left(\dfrac{1}{y^{*,n}_1},\dots,\dfrac{1}{y^{*,n}_N}\right)^{\!\mathsf{T}}\right), \\
        & \tilde{\bm{\alpha}}^{*,n} = \argmin_{\bm{a} \in A} \left\{\mathcal{I} \left[ V^n\right]((I-\Delta t \ M(\bm{y}^{*,n},\bm{a}))^{-1} \cdot \bm{y}^{*,n})  + \Delta t \,\ell(\bm{y}^{*,n},\bm{a},t^n)\right\} \\
        & \bm{y}^{*,n+1} = (I-\Delta t \ M(\bm{y}^{*,n}, \tilde{\bm{\alpha}}^{*,n}))^{-1} \cdot \bm{y}^{*,n}
    \end{aligned}
\right.
\end{equation}
where $\bm{y}^{*,0}=\bm{y}^{0}$, $\bm{y}^{*,n}=(y_1^{*,n},\dots,y_N^{*,n})\tpose\approx\bm{y}^{*}(t_n)=(y_1^{*}(t_n),\dots,y_N^{*}(t_n))\tpose$ and $\tilde{\bm{\alpha}}^{*,n}\approx \tilde{\bm{\alpha}}^*(t_n),$ for $n=0,\dots,\bar{n}-1.$

\begin{remark} \label{rmk:nointerp}
    We point out that employing more accurate Modified Patankar (MP) discretizations in place of the MPE scheme is feasible but requires specific considerations. First, since the entire reconstruction procedure relies on \eqref{eqn:MPSLscheme} and the spatial interpolation operator employed therein, increasing the order of the time integrator in \eqref{eqn:feedback_rec_MPE} may not be beneficial in terms of the overall accuracy-complexity trade-off. Furthermore, since the function $\tilde{\bm{\alpha}}^* \in \mathcal{A}$ is assumed to be just measurable, numerical methods requiring time interpolation of the optimal control may prove inadequate. As a matter of fact, interpolating over a jump could artificially regularize the control, resulting in a potential loss of optimality. For this reason, the stageless modified Patankar linear multistep methods introduced in \cite{MPLM} may be a viable alternative to other unconditionally positive and conservative MP schemes. 
\end{remark}

\subsection{Pseudo-codes}
We conclude this section with Algorithms \ref{alg: MPSL} and \ref{alg: feedback}, regarding a pseudo-code implementation of the MPSL scheme \eqref{eqn:MPSLscheme} and the reconstruction procedure \eqref{eqn:feedback_rec_MPE}, respectively. There, the parallelization in space is conceptually represented with a \textbf{parfor} instruction. As already outlined, the minimization over the compact set $A\subset \mathbb{R}^M$ is performed, for the sake of simplicity and to provide a fully discrete algorithm, by direct comparison.
\begin{algorithm}
\caption{Fully discrete MPSL algorithm} \label{alg: MPSL}
    \begin{algorithmic}[1]
        \State Compute a uniform time discretization $\{t_n\}_{n=0}^{\bar{n}}$ with step $\Delta t$ and $t_f=\bar{n} \Delta t.$
        \State Compute a space grid $\{\bm{x}_i\}_{i=0}^{\bar{I}} \subset \Omega_0$ and a control grid $\{\bm{a}_j\}_{j=0}^{\bar{J}} \subset A$
        \For {$i=0,\ldots,\bar{I}$}
            \State $V^{\bar{n}}_i \gets \varphi(\bm{x}_i)$ 
        \EndFor
        \For {$n=\bar{n},\ldots,1$} 
            \ParFor {$i=0,\ldots,\bar{I}$} 
                \State $V^* \gets +\infty$
                \For{$j=0,\ldots,\bar{J}$}
                    \State compute $M(\bm{x}_i,\bm{a}_j)$ as in \eqref{eqn:MPSLscheme}
                    \State $\bm{y} \gets$ solve $[I-\Delta t \ M(\bm{x}_i,\bm{a}_j)]\, \bm{y} = \bm{x}_i$ 
                    \State $\overline{V} \gets \mathcal{I} \left[ V^n\right](\bm{y}) + \Delta t \,\ell(\bm{x}_i,\bm{a}_j,t^n)$ 
                    \If {$\overline{V} < V^*$ } 
                        \State $V^* \gets \overline{V}$
                    \EndIf
                \EndFor
                \State $V^{n-1}_i \gets V^*$
            \EndParFor
        \EndFor
    \end{algorithmic}
\end{algorithm}
\begin{algorithm}
\caption{Optimal control and trajectory reconstruction algorithm} \label{alg: feedback}
    \begin{algorithmic}[1]
        \Require $\{t_n\}_{n=0}^{\bar{n}},$ $\{\bm{x}_i\}_{i=0}^{\bar{I}}$ and $\{V_i^n\}_{0\leq i \leq \bar{I}}^{0\leq n \leq \bar{n}},$ from Algorithm \ref{alg: MPSL}
        \Require initial data $\bm{y}^0 \in \Omega_0$
        \State Compute a control grid $\{\bm{\hat{a}}_k\}_{k=0}^{\bar{K}} \subset A$ 
        \State $\bm{y}^{*,0} \gets \bm{y}^0$
        \For {$n=0, \ldots,\bar{n}-1$}
            \State $V^* \gets +\infty$
            \For{$k=0,\ldots,\bar{K}$}
                \State compute $M(\bm{x}_i, \bm{a}_k )$ as in \eqref{eqn:MPSLscheme}
                \State $\bm{y} \gets$ solve $[I-\Delta t \ M(\bm{x}_i,\bm{a}_k )]\, \bm{y} = \bm{x}_i$
                \State $\overline{V} \gets \mathcal{I} \left[ V^n\right](\bm{y}) + \Delta t \,\ell(\bm{x}_i, \bm{a}_k ,t^n)$
                \If {$\overline{V} < V^*$ }
                    \State $\bm{y}^{*,n+1} \gets \bm{y}$ 
                    \State $\bm{a}^{*,n} \gets \bm{a}_k$ 
                \EndIf
            \EndFor
        \EndFor
    \end{algorithmic}
\end{algorithm}

\begin{remark}
In order to approximate the feedback map $(\bm{x},t) \mapsto \bm{\alpha}^*(\bm{x},t)$, one can modify Algorithm \ref{alg: MPSL} by also saving $\bm{a}^*(x_i,t^n) \gets \bm{a}_j$ right after line 14. The resulting discrete feedback could then be used in place of $V$ in Algorithm \ref{alg: feedback}, yielding an alternative yet equivalent method to compute the optimal trajectory and open-loop control.
\end{remark}

\section{Case studies and numerical simulations}\label{sec:numerical}

In this section, we present two specific CPDS with realistic applications and formulate the corresponding finite horizon optimal control problems. To numerically solve these problems, we employ a \texttt{CUDA C++} implementation of the MPSL scheme \eqref{eqn:MPSLscheme} and a \texttt{C++} implementation of the control-trajectory reconstruction \eqref{eqn:feedback_rec_MPE}, each equipped with a multilinear interpolation operator $\mathcal{I}$. For both case studies, we compare the results of the MPSL scheme with those of a classical SL approach by directly computing the discrete cost functional along the optimal trajectories provided by the two algorithms. The aim of this comparison is to provide numerical evidence that the conservative MPSL scheme produces higher-quality solutions, i.e. trajectories that are closer to the final objective at a lower cost. 

The first problem deals with chains of chemical reactions, more in particular enzyme-catalyzed reactions. Processes of this kind are crucial for the biochemical industry, because they offer high specificity and efficiency, enabling the production of complex molecules with fewer byproducts and lower energy consumption. Additionally, these reactions often occur under mild conditions, reducing the need for harsh chemicals and extreme temperatures, which enhances sustainability and safety \cite{biochem,SR18,WALSH20101}.

The second case study is about compartmental epidemic models, which provide insights on the evolution of infectious diseases and enable predictions of how they spread. In this context, control theory may represent a reliable tool to minimize the direct health impacts, such as morbidity and mortality, and to reduce the broader societal and economic disruptions. For instance, effective control measures can prevent healthcare systems from being overwhelmed (cf.  \cite{BUONOMO2014,Bolzoni2021,BOLZONI2019,Cacace_Oliviero} and references therein).

We highlight that several contributions in the literature have focused on developing positivity-preserving methods for the aforementioned applications (see, for instance, \cite{Blanes,TAKACS2024,MPV_Axioms,MPV_Mixing,Zhang2024-br,Zafar,MPVNSFD,Lubuma_Chem,Mario_Pezzella}).

\subsection{Enzyme-catalyzed biochemical reaction}

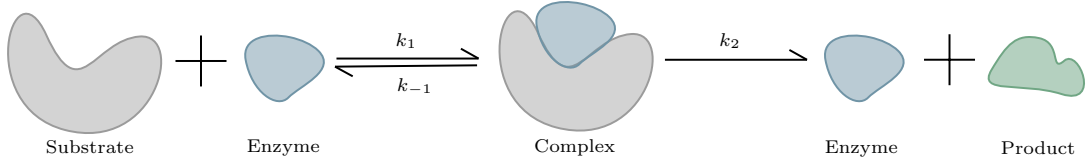
\begin{figure}
\begin{center} \hspace{-26pt}
\tikzset{every picture/.style={line width=0.75pt}} 
\begin{tikzpicture}[x=0.75pt,y=0.75pt,yscale=-1,xscale=1]
\draw  [color={rgb, 255:red, 155; green, 155; blue, 155 }  ,draw opacity=1 ][fill={rgb, 255:red, 155; green, 155; blue, 155 }  ,fill opacity=0.44 ] (72.68,30.03) .. controls (89.83,67.76) and (96.69,32.61) .. (119.84,30.89) .. controls (143,29.18) and (128.42,90.06) .. (82.11,78.05) .. controls (35.81,66.05) and (55.53,-7.7) .. (72.68,30.03) -- cycle ;
\draw   (137.89,43.24) -- (163.62,43.24)(150.75,30.38) -- (150.75,56.1) ;
\draw    (219.08,42.01) -- (290.74,42.01) ;
\draw    (281.85,37.04) -- (290.74,42.01) ;
\draw  [color={rgb, 255:red, 112; green, 148; blue, 166 }  ,draw opacity=1 ][fill={rgb, 255:red, 112; green, 148; blue, 166 }  ,fill opacity=0.48 ] (194.85,61.68) .. controls (189.13,67.54) and (176.27,59.17) .. (173.26,42.6) .. controls (171.23,31.43) and (184.06,28.79) .. (195.7,31.27) .. controls (201.32,32.47) and (206.66,34.87) .. (209.91,38.08) .. controls (219.88,47.94) and (200.57,55.81) .. (194.85,61.68) -- cycle ;
\draw    (289.99,45.43) -- (218.33,46.1) ;
\draw    (227.27,50.98) -- (218.33,46.1) ;
\draw  [color={rgb, 255:red, 155; green, 155; blue, 155 }  ,draw opacity=1 ][fill={rgb, 255:red, 155; green, 155; blue, 155 }  ,fill opacity=0.44 ] (320.55,28.78) .. controls (337.7,66.51) and (344.56,31.35) .. (367.71,29.64) .. controls (390.86,27.92) and (376.28,88.8) .. (329.98,76.8) .. controls (283.67,64.79) and (303.39,-8.95) .. (320.55,28.78) -- cycle ;
\draw  [color={rgb, 255:red, 112; green, 148; blue, 166 }  ,draw opacity=1 ][fill={rgb, 255:red, 112; green, 148; blue, 166 }  ,fill opacity=0.48 ] (341.68,44.73) .. controls (335.96,50.6) and (323.11,42.22) .. (320.1,25.66) .. controls (318.07,14.49) and (330.9,11.84) .. (342.53,14.33) .. controls (348.15,15.53) and (353.5,17.93) .. (356.74,21.14) .. controls (366.71,31) and (347.4,38.87) .. (341.68,44.73) -- cycle ;
\draw    (384.74,42.63) -- (456.4,42.63) ;
\draw    (447.51,37.66) -- (456.4,42.63) ;
\draw   (515.4,42.61) -- (541.12,42.61)(528.26,29.75) -- (528.26,55.47) ;
\draw  [color={rgb, 255:red, 112; green, 166; blue, 132 }  ,draw opacity=1 ][fill={rgb, 255:red, 112; green, 166; blue, 132 }  ,fill opacity=0.52 ] (584.79,59.94) .. controls (564.71,53.67) and (541.84,60.4) .. (549.51,47.62) .. controls (557.18,34.84) and (561.86,29.22) .. (573.5,31.71) .. controls (585.13,34.19) and (578.74,49.48) .. (586.16,43.1) .. controls (593.58,36.73) and (604.87,66.22) .. (584.79,59.94) -- cycle ;
\draw  [color={rgb, 255:red, 112; green, 148; blue, 166 }  ,draw opacity=1 ][fill={rgb, 255:red, 112; green, 148; blue, 166 }  ,fill opacity=0.48 ] (487.25,61.68) .. controls (481.53,67.54) and (468.68,59.17) .. (465.67,42.6) .. controls (463.64,31.43) and (476.47,28.79) .. (488.1,31.27) .. controls (493.72,32.47) and (499.06,34.87) .. (502.31,38.08) .. controls (512.28,47.94) and (492.97,55.81) .. (487.25,61.68) -- cycle ;
\draw (71.02,81.56) node [anchor=north west][inner sep=0.75pt]   [align=left] {{\scriptsize Substrate}};
\draw (172.33,81.15) node [anchor=north west][inner sep=0.75pt]   [align=left] {{\scriptsize Enzyme}};
\draw (317.33,81.15) node [anchor=north west][inner sep=0.75pt]   [align=left] {{\scriptsize Complex}};
\draw (463.83,81.56) node [anchor=north west][inner sep=0.75pt]   [align=left] {{\scriptsize Enzyme}};
\draw (552.72,81.98) node [anchor=north west][inner sep=0.75pt]   [align=left] {{\scriptsize Product}};
\draw (247.41,27.59) node [anchor=north west][inner sep=0.75pt]  [font=\scriptsize] [align=left] {$\displaystyle k_{1}$};

\draw (248.04,50.18) node [anchor=north west][inner sep=0.75pt]  [font=\scriptsize] [align=left] {$\displaystyle k_{-1}$};
\draw (410.56,27.59) node [anchor=north west][inner sep=0.75pt]  [font=\scriptsize] [align=left] {$\displaystyle k_{2}$};
\end{tikzpicture}
\end{center}\vspace{-0.6cm}
\caption{Visual representation of an enzymatic reaction mechanism.}
\label{fig:Chimica}
\end{figure}

Our first case study addresses  the control of a single-substrate enzyme catalyzed biochemical reaction.
In particular, our objective is to control the chain of processes schematically represented as follows
\begin{equation*}
\ce{E + S  <=>[$k_1$][$k_{-1}$]  C  ->[$k_2$]  E + P ,}
\end{equation*}
where an enzyme (E) bonds to a substrate (S), forming a complex (C) through a reversible reaction and catalysing a second reaction that yields the final product (P) and makes the enzyme available again (see Figure \ref{fig:Chimica} for a schematic representation). In this context $k_i,$ for $i=-1,1,2,$ denote the kinetic reaction rates.

The mathematical description of the aforementioned scenario is achieved through the Michaelis–Menten kinetic model (see, for instance, \cite[Section 3.1]{Ingalls2013} and \cite{Michaelis2})
\begin{equation}\label{eq:MM_original}
    \left\{
    \begin{aligned}
        s^{\prime}(t) &= - k_1 s(t)e(t) + k_{-1} c(t), \\
        e^{\prime}(t) &= - k_1 s(t) e(t) + k_{-1} c(t) + k_2  c(t), \\
        c^{\prime}(t) &= k_1 s(t) e(t) - k_{-1} c(t) -  k_2 c(t), \\
        p^{\prime}(t) &= k_2 c(t),
    \end{aligned}
    \right.
\end{equation}
where the non-negative functions $s(t),$ $e(t),$ $c(t)$ and $p(t)$ represent the relative concentrations of substrate, enzyme, complex and product, at time $t\geq t_0\geq 0,$ respectively. Here, the initial values $s_0=s(t_0),$ $e_0=e(t_0),$ $c_0=c(t_0)$ and $p_0=p(t_0)$ are considered given. 
Since the enzyme is not consumed throughout the reactions, its total concentration $e_{tot}$ remains constant and therefore
\begin{equation}\label{eq:enzima_totale}
    e(t) = e_{tot} - c(t), \qquad \mbox{for all } t \geq t_0.
\end{equation}
It turns out that the system \eqref{eq:MM_original} can be rewritten as a positive and fully conservative PDS of the form \eqref{eq:PDS_not_controlled} by taking  $\bm{y}(t)=(s(t),c(t),p(t))\tpose,$ $\bm{y}^0=(s_0,c_0,p_0)\tpose$ and
\begin{equation}\label{eq:enzyme_PDS}
    P(\bm{y}(t))=\begin{pmatrix}
        0 & k_{-1}c(t) & 0 \\
         k_1 s(t)(e_{tot} - c(t)) & 0 & 0 \\
          0 & k_2c(t) & 0 
    \end{pmatrix}=D\tpose(\bm{y}(t)).
\end{equation}
In this case it is evident that the properties \ref{ass:A1}-\ref{ass:A4} are satisfied. Furthermore, given the relation \eqref{eq:enzima_totale}, the linear invariant in \eqref{eq:conservation_law_PDS} corresponds to Lavoisier's law of conservation of mass.

Due to the dependency of the reaction dynamics on the environmental temperature, a practical approach to physically control the system \eqref{eq:enzyme_PDS} is by temperature modulation. Let $T(t)>0$ represent the absolute temperature at time $t.$ Following the arguments in \cite[Section 10.9]{Atkins-DePaula} and \cite[Page 174]{Laidler}, we adopt the modified Arrhenius law  to express the reaction rates as follows
\begin{equation}\label{eq:Arrhenius}
    k_i \sqrt{T(t)} \, \mathsf{e}^{-\frac{E^a_i}{R\, T(t)}}, \qquad \qquad i=-1,1,2.
\end{equation}
Here, the constants $k_i$ and $E_i^a,$ $i=-1,1,2,$ serve as the pre-exponential factors and the activation energies  \cite{Laidler} of the three reactions, respectively. Collectively referred to as the Arrhenius parameters, they might be estimated, for instance, from experimental data. Furthermore, $R$ denotes the universal gas constant \cite[Table 1.1]{Atkins-DePaula}. It then appears natural to consider the absolute temperature $0<T_{\min}\leq T(t)\leq T_{\max}$ as a control function and to define, starting from \eqref{eq:enzyme_PDS}, the controlled PDS 
\begin{equation}\label{eq:MM_controlled_all}
    \left\{
     \begin{aligned}
        s^{\prime}(t) &= \sqrt{T(t)} \left(- k_1  \mathsf{e}^{-\frac{E^a_1}{R T(t)}} s(t) (e_{\text{tot}} - c(t)) + k_{-1}  \mathsf{e}^{-\frac{E^a_{-1}}{R T(t)}} c(t)\right), \\
        c^{\prime}(t) &= \sqrt{T(t)} \left( k_1 \mathsf{e}^{-\frac{E^a_1}{R T(t)}} s(t) (e_{\text{tot}} - c(t)) - \left( k_{-1} \mathsf{e}^{-\frac{E^a_{-1}}{R T(t)}} + k_2 \mathsf{e}^{-\frac{E^a_2}{R T(t)}} \right) c(t)\right), \\
        p^{\prime}(t) &= k_2 \sqrt{T(t)} \ \mathsf{e}^{-\frac{E^a_2}{R T(t)}} c(t).
    \end{aligned}
    \right.
\end{equation}
For the sake of conciseness, we have opted in \eqref{eq:MM_controlled_all} not to report the explicit dependence of the state functions $s(t),$ $c(t),$ and $p(t)$ on the control $T(t).$
The system \eqref{eq:MM_controlled_all} fits the general form of a CPDS \eqref{eq:CPDS} with $A=[T_{\min}, T_{\max}],$ the production-destruction terms given by \eqref{eq:enzyme_PDS} and 
\begin{equation*}
     \mathcal{P}(T(t))=\begin{pmatrix}
        0 & \sqrt{T(t)} \ \mathsf{e}^{-\frac{E^a_{-1}}{R\, T(t)}}\phantom{z} & 0 \\
        \sqrt{T(t)} \ \mathsf{e}^{-\frac{E^a_{1}}{R\, T(t)}} & 0\phantom{z} & 0 \\
          0 & \sqrt{T(t)} \ \mathsf{e}^{-\frac{E^a_{2}}{R\, T(t)}} \phantom{z} & 0 
    \end{pmatrix}=\mathcal{D}\tpose(T(t)).
\end{equation*}
Physically, the direct intervention on the temperature modulates the reaction rates, thereby enabling the control of the overall process to attain the objective of maximizing the product in the time interval $[t_0,t_f]$. From a mathematical point of view, we introduce the space of admissible controls 
\begin{equation*}
    \mathcal{A}=\{T : [t_0, t_f] \to A=[T_{\min}, T_{\max}] \; |\; T \text{ measurable} \}
\end{equation*}
and define, for $T\in \mathcal{A},$ the running and final costs as follows
\begin{equation*}
    \ell(T(t))= w_{1}\, \left( \frac{T(t)-T_{amb}}{T_{\max}} \right)^2, \qquad \varphi(p(t_f))= w_{2}\, (1-p(t_f))^2,
\end{equation*}
where $T_{min}=263.15$ K, $T_{amb}=293.15$ K and $T_{max}=373.15$ K represent the minimum, ambient and maximum temperatures, respectively. Here, the weights are set to $w_1=\frac12$ and $w_2=20$ to prioritize the maximization of the final product over regularizing the temperature.
Hence, our ODE-constrained, finite horizon, optimal control problem reads
\begin{equation}\label{eq:Opt_Cntr_ENZIMA}
    \begin{aligned}
    & \underset{T\in \mathcal{A}}{\text{minimize}}
    & & \int_{t_0}^{t_f} \frac{1}{2}\, \left( \frac{T(t)-T_{amb}}{T_{\max}} \right)^2  dt \ + 20\, (1-p(t_f;T(t_f)))^2, \\
    & \text{subject to} 
    & & \eqref{eq:MM_controlled_all}.
    \end{aligned}
\end{equation}
\subsubsection{Numerical Simulations}
We report some simulation outcomes for the optimal control problem \eqref{eq:Opt_Cntr_ENZIMA}. Specifically, we consider the CPDS \eqref{eq:MM_controlled_all} with $\bm{y}^0=(0.7,0,0)\tpose$ and the parameters detailed in Table \ref{tab:enzyme}, that are meant to be realistic but are not inferred from experimental data.
\begin{table}[H]
    \centering
    \begin{tabular}{|c|c|c|c|c|c|c|c|c|c|}
        \hline
        \textit{Parameter} & $t_0$ & $t_f$ & $k_1$ & $k_{-1}$ & $k_2$ & $E_1^a / R$ & $E_{-1}^a / R$ & $E_2^a / R$ & $e_{\text{tot}}$ \\
        \hline
        \textit{Value} & 0 & 30 & 0.04 & 0.03 & 0.035 & 200 & 220 & 190 & 0.3 \\
        \hline
    \end{tabular}
    \caption{Parameter values used for the numerical simulation and control of the CPDS \eqref{eq:MM_controlled_all}.}
    \label{tab:enzyme}
\end{table}

We discretize the cube $[0,1]^3\supset \Omega_0$ with a uniform grid with step $\Delta x=1/750,$ i.e. made of $751^3\approx 4.24\cdot 10^8$ nodes. Then, the solution of the HJB equation \eqref{eq:HJB} is approximated by \eqref{eqn:MPSLscheme} with $\Delta t = 30/100,$ while for the optimal control and trajectory reconstruction the procedure \eqref{eqn:feedback_rec_MPE} is employed. Furthermore, the control set $A$ is uniformly discretized with $10^3$ points. We emphasize that the choice $\Delta t > \Delta x,$ which could be unfeasible in a finite-difference setting due to CFL-like stability conditions, is here adopted to reduce numerical dissipation effects.
 
\begin{figure}[H]
    \centering
    \includegraphics[width=1\linewidth]{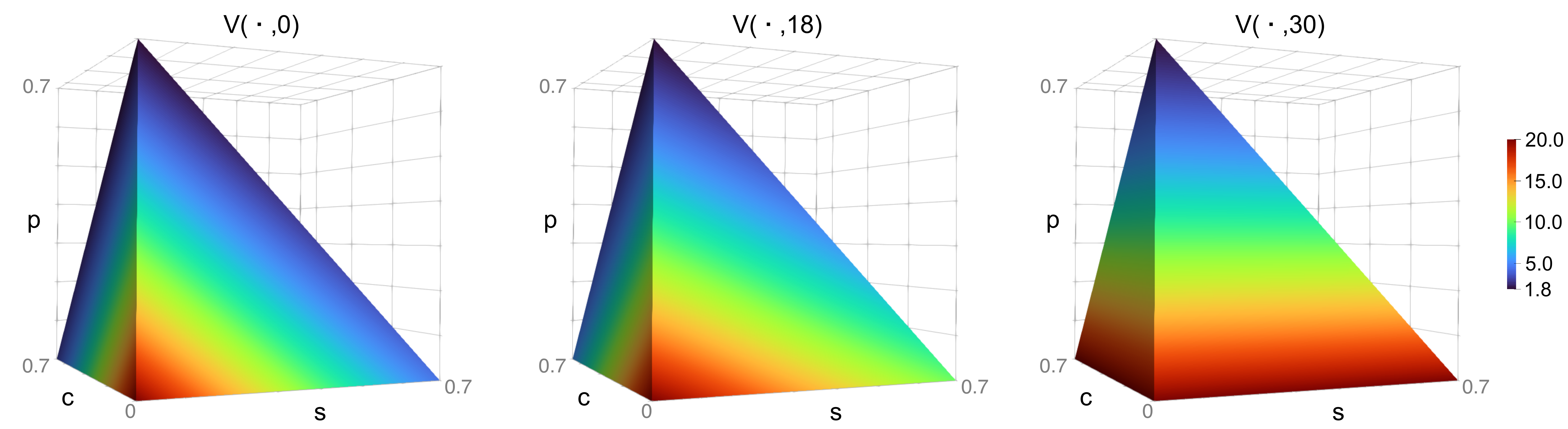}
    \caption{Discrete value function in the positively invariant region $\Omega_0=\{\bm{x}\in\mathbb{R}^3 \ : \ \bm{0}\leq \bm{x} \leq \bm{e}\tpose \bm{y}^0  \}$ for $t=0$ (left), $t=18$ (center) and $t=30$ (right).}
    \label{fig:value_time_chimica}
\end{figure}

\begin{figure}[H]
    \centering
    \includegraphics[width=0.33\linewidth]{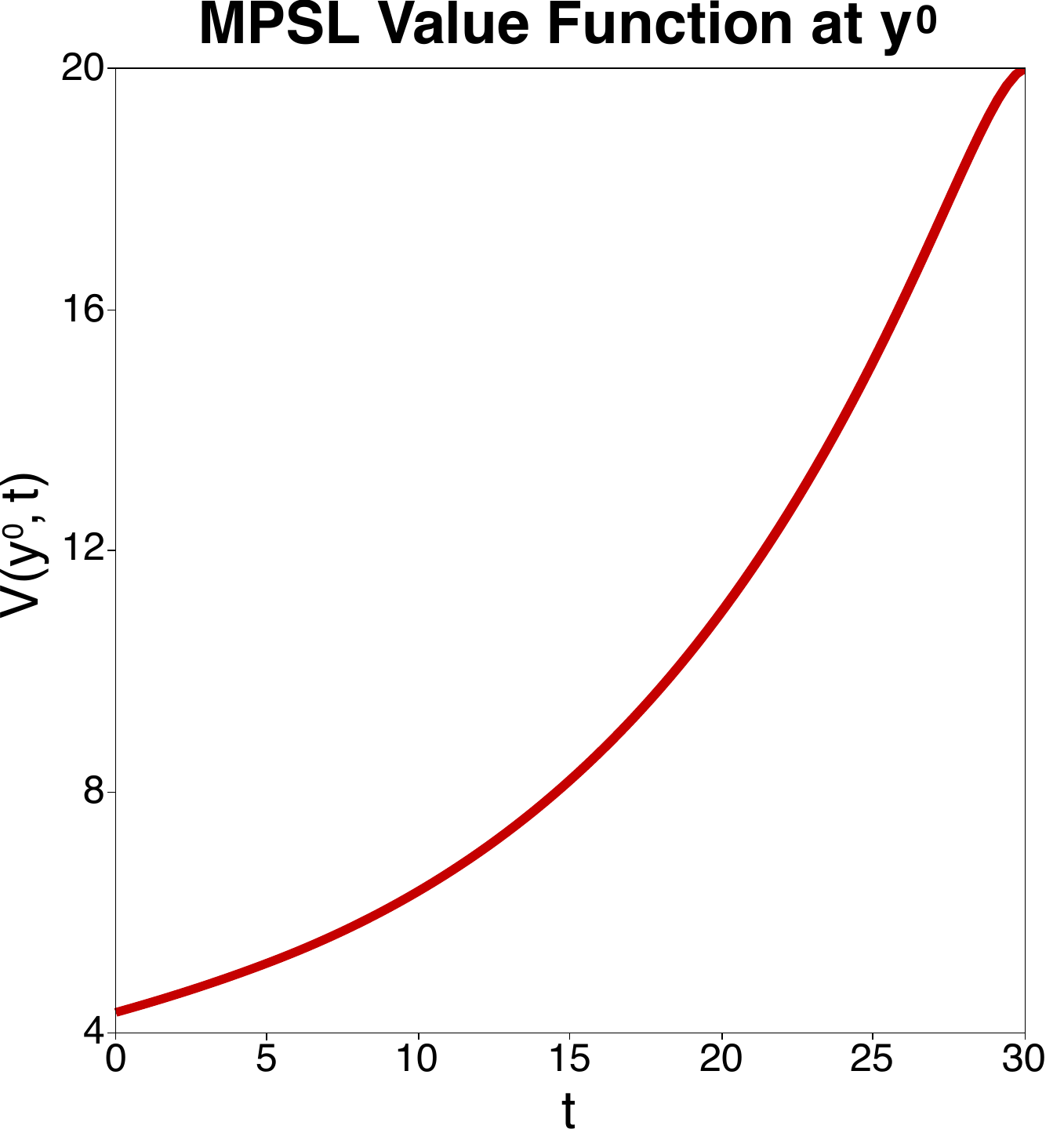}
    \caption{Time-evolution of the discrete value function at the initial data $\bm{y}^0=(0.7, 0, 0)\tpose$.} 
    \label{fig:value_chimica_y0}
\end{figure}

The discrete value function obtained with the MPSL scheme \eqref{eqn:MPSLscheme} is reported in Figure \ref{fig:value_time_chimica} for three different values of $t.$ Moreover, we show in Figure \ref{fig:value_chimica_y0} the time-evolution of the approximated value function at $\bm{x}= \bm{y}^0 = (0.7,0,0)\tpose$. It is worth to note that in our finite horizon setting, $V(\bm{y}^0,t)$ is a non-decreasing function of time, since, as $t$ increases, the effort to achieve the objective must grow, resulting in a higher optimal cost. We remark that the value function evolves backward in time, therefore the plots of Figures \ref{fig:value_time_chimica} and \ref{fig:value_chimica_y0} are meant to be read from right to left.

\begin{figure}[H]
    \centering
    \includegraphics[width=0.7\linewidth]{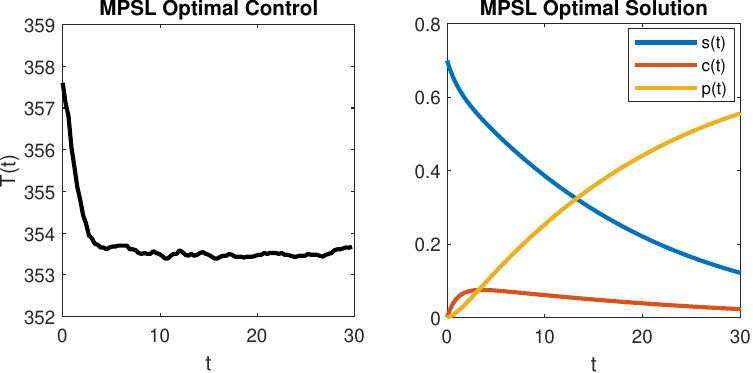}
    \caption{Numerical solution of the optimal control problem \eqref{eq:Opt_Cntr_ENZIMA} computed by the MPSL scheme \eqref{eqn:MPSLscheme} with the value and feedback reconstruction \eqref{eqn:feedback_rec_MPE}.}
    \label{fig:mpsl_chimica}
\end{figure}

\begin{figure}[H]
    \centering
    \includegraphics[width=0.7\linewidth]{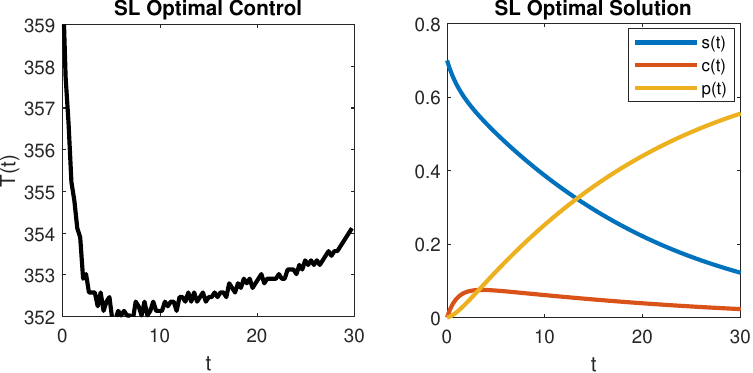}
    \caption{Numerical solution of the optimal control problem \eqref{eq:Opt_Cntr_ENZIMA} computed by the classical SL scheme with the value and feedback reconstruction \eqref{eqn:feedback_rec_MPE}.}
    \label{fig:sl_chimica}
\end{figure}

As for the optimal controls and trajectories, those derived with the MPSL scheme \eqref{eqn:MPSLscheme} and the reconstruction procedure \eqref{eqn:feedback_rec_MPE} are reported in Figure \ref{fig:mpsl_chimica}, whereas those given by a classical SL method, employing an explicit Euler integrator to trace the characteristics, are presented in Figure \ref{fig:sl_chimica}. First of all, as shown in Table \ref{tab:costi_chimica}, both approaches for solving \eqref{eq:Opt_Cntr_ENZIMA} provide better solutions compared to the base case $T(t)\equiv T_{amb},$ leading to higher values for the final product and lower optimal trajectory costs. Furthermore, when compared to the classical SL method, we notice that the two optimal control functions $T(t)$ are not identical. Although both end at around the same temperature, the one in Figure \ref{fig:mpsl_chimica} is flatter and more regular. This difference is not only qualitative. In fact, the advantages of the MPSL scheme are twofold: it results in a more substantial increase in the value of $p(t_f)$ and achieves a more considerable reduction in the corresponding objective functional $J_{\bm{y}^0,t_0} (T).$

\begin{table}[H]
    \centering    
\begin{tabular}{|c|c|c|c|c|}
    \multicolumn{5}{c}{Base Case} \\
    \cline{1-5}
      & \multicolumn{2}{c|}{$p(t_f)=0.500824$} & \multicolumn{2}{c|}{$J_{\bm{y}^0,t_0} (T)=4.983513$ } \\ 
    \cline{1-5} 
    \multicolumn{5}{c}{\rule{0pt}{0.6cm}Optimal Control Setting} \\
     \hline
    \textit{Scheme} & $p(t_f)$ & $p_{\%}$ & $J_{\bm{y}^0,t_0} (T)$ & $J_{\%}$ \\
    \hline
    MPSL & 0.556064 & $+11.03\%$ & 4.337023 & $-12.97\%$\\
    \hline
    SL & 0.555476 & $+10.91\%$  & 4.337484 & $-12.96\%$\\
    \hline
\end{tabular}
\caption{Final product and cost functional evaluated along the optimal trajectories for system \eqref{eq:MM_controlled_all} in both the base case ($T(t) \equiv T_{amb}$) and the controlled scenarios. Here, $p_{\%}$ and $J_{\%}$ denote the percentage variation of $p(t_f)$ and $J_{\bm{y}^0,t_0}$ relative to the base case, respectively.}
\label{tab:costi_chimica}
\end{table}

\subsection{Infectious diseases spreading}
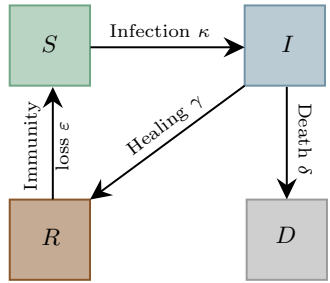
\begin{figure}[H]
\begin{center} \hspace{-25pt}
\tikzset{every picture/.style={line width=0.75pt}} 
\begin{tikzpicture}[x=0.75pt,y=0.75pt,yscale=-1,xscale=1]
\draw  [color={rgb, 255:red, 112; green, 148; blue, 166 }  ,draw opacity=1 ][fill={rgb, 255:red, 112; green, 148; blue, 166 }  ,fill opacity=0.48 ] (198.5,52.5) -- (239,52.5) -- (239,93) -- (198.5,93) -- cycle ;
\draw  [color={rgb, 255:red, 112; green, 166; blue, 132 }  ,draw opacity=1 ][fill={rgb, 255:red, 112; green, 166; blue, 132 }  ,fill opacity=0.48 ] (80,52.5) -- (120.5,52.5) -- (120.5,93) -- (80,93) -- cycle ;
\draw  [color={rgb, 255:red, 155; green, 155; blue, 155 }  ,draw opacity=1 ][fill={rgb, 255:red, 155; green, 155; blue, 155 }  ,fill opacity=0.48 ] (199.5,150.5) -- (240,150.5) -- (240,191) -- (199.5,191) -- cycle ;
\draw  [color={rgb, 255:red, 139; green, 87; blue, 42 }  ,draw opacity=1 ][fill={rgb, 255:red, 139; green, 87; blue, 42 }  ,fill opacity=0.5 ] (80,150.5) -- (120.5,150.5) -- (120.5,191) -- (80,191) -- cycle ;
\draw    (120.74,73.63) -- (194.5,73.63) ;
\draw [shift={(198.5,73.63)}, rotate = 180] [fill={rgb, 255:red, 0; green, 0; blue, 0 }  ][line width=0.08]  [draw opacity=0] (10.72,-5.15) -- (0,0) -- (10.72,5.15) -- (7.12,0) -- cycle    ;
\draw    (219.74,93.63) -- (219.74,146.4) ;
\draw [shift={(219.74,150.4)}, rotate = 270] [fill={rgb, 255:red, 0; green, 0; blue, 0 }  ][line width=0.08]  [draw opacity=0] (10.72,-5.15) -- (0,0) -- (10.72,5.15) -- (7.12,0) -- cycle    ;
\draw    (101.74,150.4) -- (101.74,96.63) ;
\draw [shift={(101.74,93.63)}, rotate = 90] [fill={rgb, 255:red, 0; green, 0; blue, 0 }  ][line width=0.08]  [draw opacity=0] (10.72,-5.15) -- (0,0) -- (10.72,5.15) -- (7.12,0) -- cycle    ;
\draw    (122.91,148.72) -- (198.5,93) ;
\draw [shift={(120.5,150.5)}, rotate = 323.6] [fill={rgb, 255:red, 0; green, 0; blue, 0 }  ][line width=0.08]  [draw opacity=0] (10.72,-5.15) -- (0,0) -- (10.72,5.15) -- (7.12,0) -- cycle    ;
\draw (94.5,66.5) node [anchor=north west][inner sep=0.75pt]  [font=\small] [align=left] {$\displaystyle S$};
\draw (215.5,66.5) node [anchor=north west][inner sep=0.75pt]  [font=\small] [align=left] {$\displaystyle I$};
\draw (212.5,162.5) node [anchor=north west][inner sep=0.75pt]  [font=\small] [align=left] {$\displaystyle D$};
\draw (94.5,163.5) node [anchor=north west][inner sep=0.75pt]  [font=\small] [align=left] {$\displaystyle R$};
\draw (85.83,147.50) node [anchor=north west][inner sep=0.75pt]  [rotate=-270] [align=center] {{\scriptsize Immunity }\\{\scriptsize loss $\varepsilon$}};
\draw (135.97,123.04) node [anchor=north west][inner sep=0.75pt]  [rotate=-322.6] [align=left] {{\scriptsize Healing $\gamma$}};
\draw (128.83,60.65) node [anchor=north west][inner sep=0.75pt] [rotate=0]  [align=left] {{\scriptsize Infection $\kappa$}};
\draw (233.83,100.00) node [anchor=north west][inner sep=0.75pt]  [rotate=-90] [align=left] {{\scriptsize Death $\delta$}\\};
\end{tikzpicture}
\end{center}\vspace{-0.6cm}
\caption{Visual representation of the infection pathway.}
\label{fig:SIRD}
\end{figure}

In our second case study, we explore the possibility to control the spread of infectious diseases through the implementation of restrictive measures. We consider a closed population for which demographic turnover is neglected and individuals are categorized into four distinct compartments: susceptible ($S$, those at risk of infection), infective ($I,$ those who have contracted and can spread the disease), recovered ($R,$ those who have overcome the infection) and deceased ($D$). In this scenario, we assume that the immunity following the recovery is temporary, allowing for the possibility of reinfection. Additionally, we account for psychological effects, letting the population's behavior vary in response to changes in the number of infections. We refer to Figure \ref{fig:SIRD}  for a schematic representation of the infection pathway.

Let $s(t)$, $i(t)$, $r(t)$ and $d(t)$ represent the fraction of the total population that belongs, at time $t\geq t_0,$ to the compartment $S$, $I$, $R$ and $D,$ respectively. Given the initial values $s_0=s(t_0),$ $i_0=i(t_0),$ $r_0=r(t_0)$ and $d_0=d(t_0),$ the evolution of the disease is modeled as follows
\begin{equation} \label{eq:SIRD}
    \begin{cases}
        s'(t) = - g(i(t)) s(t) + \varepsilon r(t), \\
        i'(t) =  g(i(t)) s(t) - (\gamma + \delta) i(t), \\
        r' (t) = \gamma i(t) - \varepsilon r(t), \\
        d'(t) = \delta i(t),
    \end{cases}
\end{equation}
where $\gamma,\delta, \varepsilon \in \mathbb{R}^+$ are the recovery, death and immunity loss rates. Here, the force of infection, i.e. the probability per unit of time for a susceptible individual to become infected \cite{DiekmannHeesterbeek2000}, reads 
\begin{equation*}
    g(i)=\kappa \, i^\varrho  (1-i)^\sigma, \quad \text{with} \quad \varrho\geq 1, \quad \sigma \geq 0 \quad \text{and} \quad \kappa >0.
\end{equation*}
As outlined in \cite{DOnofrio, CAPASSO197843, vandenDriessche2000} and references therein, this non-standard form for $g$ takes into account the reduction in the frequency of contacts, driven by the fear of contagion, when the proportion of infected individuals exceeds the threshold value $\frac{\varrho}{\varrho+\sigma}$. Thus, when $\varrho=1$ and $\sigma=0$, the product $g(i(t)) s(t)$ coincides with the widely used mass action incidence term. Denote by $\bm{y}(t)=(s(t),i(t),r(t),d(t))\tpose,$ for $t\geq t_0$ and by $\bm{y}^0=(s_0,i_0,r_0,d_0)\tpose.$ The differential system \eqref{eq:SIRD} then corresponds to the positive and fully conservative PDS \eqref{eq:PDS_not_controlled} with
\begin{equation}\label{eq:SIRD_PDS}
    P(\bm{y}(t))=\begin{pmatrix}
        0 & 0 & \varepsilon r(t) & 0 \\
         g(i(t)) s(t) & 0 & 0 & 0 \\
          0 & \gamma i(t) & 0 & 0 \\
          0 & \delta i(t) & 0 & 0
    \end{pmatrix}=D\tpose(\bm{y}(t)).
\end{equation}
In this context, properties \ref{ass:A1}-\ref{ass:A4} are satisfied and the linear invariant in \eqref{eq:conservation_law_PDS}  reflects the conservation of the total population size. 

Non-pharmaceutical interventions, including mobility restrictions and social distancing, limit interactions and therefore represent a direct and immediate approach to mitigate the spread of infectious diseases. We incorporate such measures into the model \eqref{eq:SIRD} by introducing a control function $u: \mathbb{R}^+_0 \to [0,1],$ which quantifies the intensity of the interventions at time $t\geq t_0.$ Specifically, we consider the controlled system
\begin{equation}\label{eq:cont_SIRD}
    \begin{cases}
        s'(t) = (u(t)-1) \, g(i(t)) s(t) + \varepsilon r(t), \\
        i'(t) = (1-u(t)) \, g(i(t)) s(t) - (\gamma + \delta) i(t), \\
        r' (t) = \gamma i(t) - \varepsilon r(t), \\
        d'(t) = \delta i(t),
    \end{cases}
\end{equation}
which follows the general CPDS formulation \eqref{eq:CPDS} with $A=[0, 1],$ the production-destruction terms defined in \eqref{eq:SIRD_PDS} and control policies $\mathcal{P}(u(t))=\mathcal{D}\tpose(u(t))\in \R^{4\times4},$ that exhibit all zero components except for
\begin{equation*}
     \mathcal{P}_{21}(u(t))=1-u(t)=\mathcal{D}_{12}(u(t)).
\end{equation*}
Here, a control strategy is designed with the objective of minimizing both the cumulative number of infections caused by the disease over the interval $[t_0,t_f]$ and the number of dead people at the final time.  Furthermore, as interventions come with inherent costs that may become too demanding, restrictions should be applied only when truly necessary.
These considerations lead to the following running and final costs
\begin{equation*}
    \ell(u(t),i(t)) = \int_{t_0}^{t_f} w_1 u^2(t) + w_2 i(t) \ dt , \qquad \quad \varphi(d(t_f)) = w_3 d(t_f),
\end{equation*}
where we choose $w_3>>w_2>>w_1$ to differently weigh the single contributions of the aforementioned cost components. Given the space of admissible controls 
\begin{equation*}
    \mathcal{A}=\{u : [t_0, t_f] \to A=[0, 1] \; |\; u \text{ measurable} \},
\end{equation*}
we then define the following ODE-constrained, finite horizon, optimal control problem
\begin{equation}\label{eq:Opt_Cntr_SIRD}
    \begin{aligned}
    & \underset{u\in \mathcal{A}}{\text{minimize}}
    & & \int_{t_0}^{t_f} 10^{-3} u^2(t) + i(t) \ dt \ +  10^4 d(t_f), \\
    & \text{subject to} 
    & & \eqref{eq:cont_SIRD}.
    \end{aligned}
\end{equation}
\subsubsection{Numerical Simulations}

Some numerical experiments addressing the optimal control problem \eqref{eq:Opt_Cntr_SIRD} are here presented. The specific parameter values for \eqref{eq:cont_SIRD} are listed in Table \ref{tab:SIRD}. Moreover, to represent a realistic situation in which restrictive measures are introduced not at the onset of the first infection but following a certain delay, we set $\bm{y}^0=(0.985, 0.007, 0.006, 0.002)^\top$, which reflects the epidemiological scenario expected after thirty days without intervention.

\begin{table}[H]
    \centering
    \begin{tabular}{|c|c|c|c|c|c|c|c|c|}
        \hline
        \textit{Parameter} & $t_0$ & $t_f$ & $\kappa$ & $\varrho$ & $\sigma$ & $\gamma$ & $\delta $ & $\varepsilon$  \\
        \hline
        \textit{Value} & 0 & 90 & 0.32 & 1 & 0.5 & 0.12 & 0.0294 & 0.0094  \\
        \hline
        \textit{Reference} & - & - & \cite{Siettos_Sird} & \cite{Liu1986} & \cite{Mickens_root} & \cite{FERNANDEZVILLAVERDE2022104318} & \cite{Siettos_Sird} & \cite{Reinfection}  \\
        \hline
    \end{tabular}
    \caption{Parameter values used for the numerical simulation and control of the CPDS \eqref{eq:cont_SIRD}.}
    \label{tab:SIRD}
\end{table} 

\begin{figure}
\centering
\includegraphics[width=0.4\linewidth]{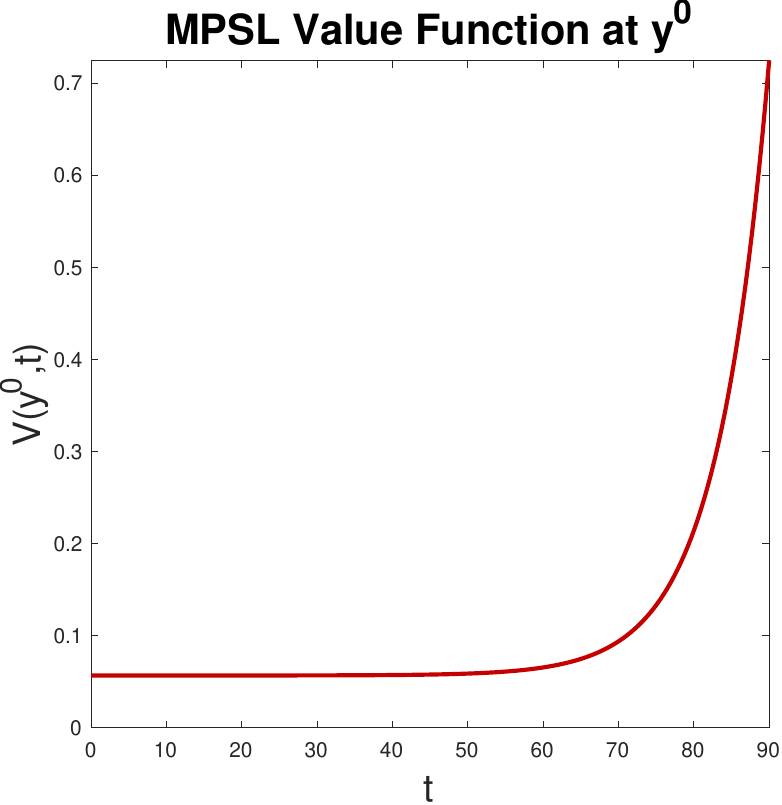}
    \caption{Time-evolution of the discrete value function at the initial data $\bm{y}^0=(0.985, 0.007, 0.006, 0.002)\tpose$.}
    \label{fig:value_sird_y0}
\end{figure}

The MPSL scheme \eqref{eqn:MPSLscheme} is adopted with $\Delta t=90/200$ and $\Delta x=1/100,$ i.e. $101^4 \approx 1.04 \cdot 10^8$ nodes in $[0,1]^4 \supset \Omega_0$, to approximate the viscosity solution of the HJB equation \eqref{eq:HJB}. Then, following a uniform discretization of the control set $A$ with a step size $\Delta a = 1/1000$, we reconstruct the optimal control and trajectory using the algorithm provided in \eqref{eqn:feedback_rec_MPE}. Figure \ref{fig:value_sird_y0} exhibits the approximated value function, computed backward in time, for the chosen initial point.

\begin{figure}[H]
    \centering
    \includegraphics[width=0.7\linewidth]{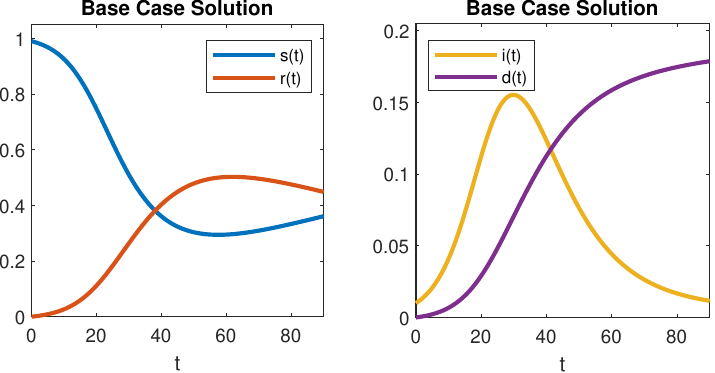}
    \caption{Numerical solution of the base case, i.e. \eqref{eq:cont_SIRD} with $u(t)\equiv 0$.}
    \label{fig:SIRD_Base_Case}
\end{figure}

\begin{figure}[H]
    \centering
    \includegraphics[width=0.95\linewidth]{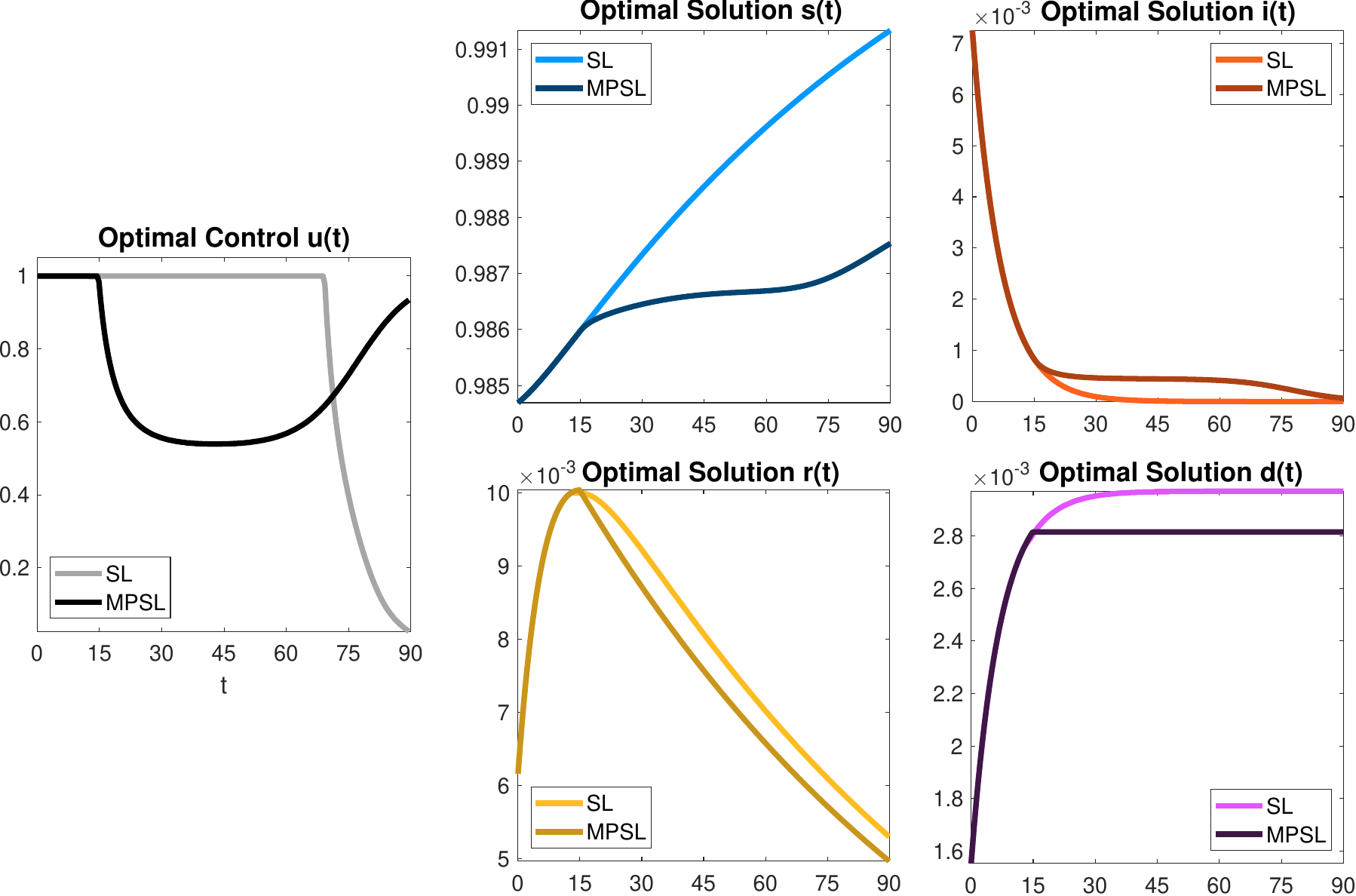}
    \caption{Numerical solution of the optimal control problem \eqref{eq:Opt_Cntr_SIRD} computed by the classical SL and the MPSL scheme \eqref{eqn:MPSLscheme} with the value and feedback}
    \label{fig:SIRD_Traj_Contr}
\end{figure}

The simulation outcomes using both \eqref{eqn:MPSLscheme} and the classical SL scheme, depicted in Figure \ref{fig:SIRD_Traj_Contr}, are compared to those of the base case (no intervention policy, i.e., $u(t)=0$ for all $t \geq t_0$), shown in Figure \ref{fig:SIRD_Base_Case}. This comparison underscores the qualitative effectiveness of the epidemic control strategy and highlights a significant reduction in the proportion of infected and deceased individuals. From a quantitative standpoint, as shown by the data in Table \ref{tab:costi_SIRD}, the MPSL method proves more accurate than the classical SL scheme, achieving a greater reduction in both the cost functional and the value of $d(t_f)$. Additionally, the control strategy suggested by \eqref{eqn:MPSLscheme} appears more realistic, advocating for the implementation of more dynamic but less strict restrictive measures, with a total lockdown period that is two months shorter than that indicated by the classical SL approach.

\begin{table}[H]
    \centering    
\begin{tabular}{|c|c|c|c|c|}
    \multicolumn{5}{c}{Base Case} \\
    \cline{1-5}
      & \multicolumn{2}{c|}{$d(t_f)=0.178829$} & \multicolumn{2}{c|}{$J_{\bm{y}^0,t_0} (u)=325.8751$ } \\ 
    \cline{1-5} 
    \multicolumn{5}{c}{\rule{0pt}{0.6cm}Optimal Control Setting} \\
     \hline
    \textit{Scheme} & $d(t_f)$ & $d_{\%}$ & $J_{\bm{y}^0,t_0} (T)$ & $J_{\%}$ \\
    \hline
    MPSL & 0.002815 & $-98.43\%$ & 0.199714 & $-99.94\%$\\
    \hline
    SL & 0.002970 & $-98.34\%$  & 0.210538 & $-99.93\%$\\
    \hline
\end{tabular}
\caption{Final deceased portion and cost functional evaluated along the optimal trajectories for system \eqref{eq:Opt_Cntr_SIRD} in both the base case ($u(t) \equiv 0$) and the controlled scenarios. Here, $d_{\%}$ and $J_{\%}$ denote the percentage variation of $d(t_f)$ and $J_{\bm{y}^0,t_0}$ relative to the base case, respectively.}
\label{tab:costi_SIRD}
\end{table}

\begin{figure}[H]
\centering
\includegraphics[width=0.4\linewidth]{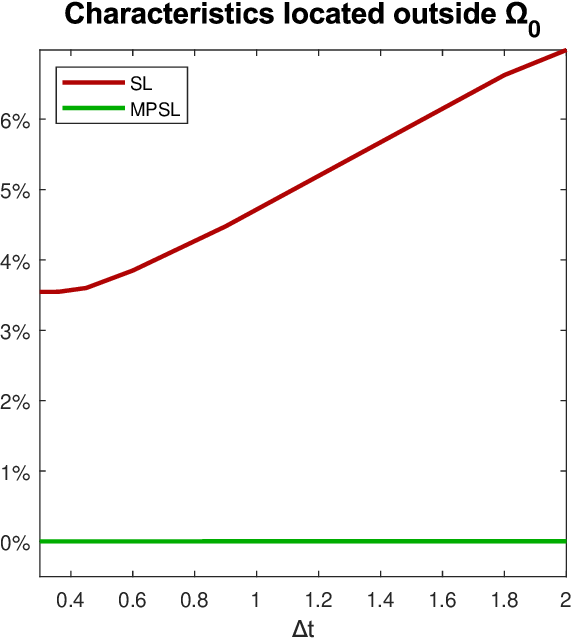}
    \caption{Percentage of characteristic curves that, starting from a node in $\omega_{\Delta x},$ are wrongly located outside $\Omega_0$ by the classical SL scheme employing forward Euler.}
    \label{fig:Uscita}
\end{figure}

We aim to provide an explanation for the inherent causes of the differences observed between the SL and MPSL control strategies illustrated in Figure \ref{fig:SIRD_Traj_Contr}. To this end, we define the subset of the computational domain
\begin{equation*}
    \omega_{\Delta x}=\{\bm{x}\in\{\bm{x}_i\}_{0\leq i \leq \bar{I}}\ |\ \dist(\bm{x},\partial\Omega_0) \leq 10 \Delta x\}
\end{equation*}
and compute the overall percentage of characteristic curves that, starting from a node in $\omega_{\Delta x},$ are wrongly located outside the positively invariant region $\Omega_0$, defined in \eqref{eq:Omega0}, after one time step. We conducted this analysis for various values of $\Delta t$ and the corresponding results are presented in Figure \ref{fig:Uscita}. As expected, by virtue of the positivity and conservativity properties inherited from the modified Patankar integrator, the MPSL scheme successfully tracks all the characteristics regardless of the time step-size. Conversely, the classical SL scheme fails to preserve the system's inherent properties due to Euler's method. As a matter of fact, due to the incorrect localization of the characteristics, certain control strategies, which would perform better if their effects were accurately captured by the discrete dynamics, are considered unfeasible. These considerations place emphasis on the limitations of classical SL schemes in the context of PDS, while also demonstrating the robustness of the dynamic programming approach. In fact, the feedback mechanism effectively compensates for the inaccuracy of the ODE solver, yielding solutions that, although suboptimal, remain acceptable.

We then conclude that the solution provided by the MPSL scheme most closely approximates the optimum.

\section{Conclusions}\label{sec:Conclusions}
In this work, we established a rigorous theoretical and numerical framework for the optimal control of positive and fully conservative Production-Destruction Systems (PDS). As a first step, we formulated the general model of a Controlled Production-Destruction System (CPDS) and investigated the requirements that control functions and policies have to satisfy in order to preserve the inherent properties of the original PDS. We then introduced a finite-horizon optimal control problem, which we addressed by means of the dynamic programming approach and the definition of a Hamilton-Jacobi-Bellman (HJB) partial differential equation. To accurately approximate the value function, i.e. the viscosity solution to the HJB equation, we proposed a parallel-in-space Modified Patankar Semi-Lagrangian (MPSL) scheme. Furthermore, we provided a straightforward positive and conservative algorithm for reconstructing the optimal control and trajectory. The MPSL procedure was implemented in \texttt{CUDA C++} and applied to two realistic case studies (enzyme-catalyzed biochemical reaction model and epidemic model), for which the mathematical properties of the continuous CPDS have a well-established physical meaning.  The numerical simulations provided evidence of the effectiveness and robustness of the proposed procedure. In particular, the comparison with classical semi-Lagrangian schemes revealed that the MPSL strategy results in a greater reduction in the objective functional. Moreover, the conservative nature of the MPSL scheme ensures the correct localization of characteristics within the CPDS's positively invariant region, leading to more realistic and accurate control strategies.

Potential advancements of the optimal control framework we presented in this paper may primarily move towards the development of high-order MPSL schemes.  In this context, more accurate modified Patankar methods that do not require time-interpolation would represent reliable integrators to correctly track the characteristic curves, in compliance with our objectives. However, in order to achieve an overall high-order MPSL procedure, they need to be complemented with a state-space reconstruction of comparable order -- such as $\mathbb{P}^d$ or $\mathbb{Q}^d$ interpolations with $d>1$ or Essentially Non-Oscillatory (ENO) reconstructions -- which may not be monotone. Therefore, proving the convergence of the resulting scheme would be more challenging, as shown in \cite{CFPS24}, where a CWENO reconstruction is employed to recover the value function on the whole computational domain. Further developments may focus on the implementation of the MPSL scheme on unstructured grids, with the aim of achieving a more accurate discretization of the positively invariant region and a more efficient use of memory.

\section*{Acknowledgements}
This work has been performed in the auspices of the \textit{Italian National Group for Scientific Computing} (GNCS) of the National Institute for Advanced Mathematics (INdAM). M. Pezzella is supported by the Project PE 0000020 CHANGES - CUP\!\_B53C22003890006, NRP Mission 4 Component 2 Investment 1.3, Funded by the European Union - NextGenerationEU. A. Oliviero is supported by the Italian PNRR fund within the doctoral project \textit{Modelli matematici per la simulazione e controllo
delle pandemie} at Sapienza University of Rome. S. Cacace is supported by the PNRR-MUR project \textit{Italian Research Center on High Performance Computing, Big Data and Quantum Computing}. A. Oliviero and S. Cacace are supported by the INdAM-GNCS Project code CUP\!\_E53C23001670001.

%
%

\bibliographystyle{plainurl}
\bibliography{Pezzella_Cacace_Oliviero}   

\begin{thebibliography}{10}

\bibitem{Siettos_Sird}
C.~Anastassopoulou, L.~Russo, A.~Tsakris, and C.~Siettos.
\newblock {Data-based analysis, modelling and forecasting of the COVID-19 outbreak}.
\newblock {\em PLOS ONE}, 15(3):1--21, 03 2020.
\newblock \href {https://doi.org/10.1371/journal.pone.0230405} {\path{doi:10.1371/journal.pone.0230405}}.

\bibitem{Lubuma_Chem}
R.~Anguelov, Y.~Dumont, and J.~M.-S. Lubuma.
\newblock {On nonstandard finite difference schemes in biosciences}.
\newblock {\em AIP Conference Proceedings}, 1487(1):212--223, 10 2012.
\newblock \href {https://doi.org/10.1063/1.4758961} {\path{doi:10.1063/1.4758961}}.

\bibitem{biochem}
V.~L. Arcus and A.~J. Mulholland.
\newblock Temperature, dynamics, and enzyme-catalyzed reaction rates.
\newblock {\em Annual Review of Biophysics}, 49(Volume 49, 2020):163--180, 2020.
\newblock \href {https://doi.org/10.1146/annurev-biophys-121219-081520} {\path{doi:10.1146/annurev-biophys-121219-081520}}.

\bibitem{Atkins-DePaula}
P.~Atkins and J.~De~Paula.
\newblock {\em Elements of physical chemistry}.
\newblock Oxford University Press, USA, 2013.

\bibitem{AKK21}
B.~Azmi, D.~Kalise, and K.~Kunisch.
\newblock Optimal feedback law recovery by gradient-augmented sparse polynomial regression.
\newblock {\em Journal of Machine Learning Research}, 22(48):1--32, 2021.

\bibitem{TAKACS2024}
{B. Máté Takács}, G.~{S. Sebestyén}, and {I. Faragó}.
\newblock High-order reliable numerical methods for epidemic models with non-constant recruitment rate.
\newblock {\em Applied Numerical Mathematics}, 206:75--93, 2024.
\newblock \href {https://doi.org/10.1016/j.apnum.2024.08.008} {\path{doi:10.1016/j.apnum.2024.08.008}}.

\bibitem{BardiCapuzzoDocetta2008}
M.~Bardi and I.~Capuzzo-Dolcetta.
\newblock {\em Optimal Control and Viscosity Solutions of {Hamilton}--{Jacobi}--{Bellman} Equations}.
\newblock Modern Birkh{\"a}user Classics. Birkh{\"a}user Boston, 2008.

\bibitem{Bellman1957}
R.~Bellman.
\newblock {\em {Dynamic Programming}}.
\newblock Princeton University Press, Princeton, NJ, 1957.

\bibitem{Blanes}
{Blanes, S.}, {Iserles, A.}, and {Macnamara, S.}
\newblock Positivity-preserving methods for ordinary differential equations.
\newblock {\em ESAIM: M2AN}, 56(6):1843--1870, 2022.
\newblock \href {https://doi.org/10.1051/m2an/2022042} {\path{doi:10.1051/m2an/2022042}}.

\bibitem{BGP56}
V.~G. Boltyanskii, R.~V. Gamkrelidze, and L.~S. Pontryagin.
\newblock On the theory of optimal processes.
\newblock In {\em Dokl. Akad. Nauk SSSR}, volume 110, pages 7--10, 1956.

\bibitem{BOLZONI2019}
L.~Bolzoni, E.~Bonacini, R.~{Della Marca}, and M.~Groppi.
\newblock Optimal control of epidemic size and duration with limited resources.
\newblock {\em Mathematical Biosciences}, 315:108232, 2019.
\newblock \href {https://doi.org/10.1016/j.mbs.2019.108232} {\path{doi:10.1016/j.mbs.2019.108232}}.

\bibitem{Bolzoni2021}
L.~Bolzoni, R.~Della~Marca, and M.~Groppi.
\newblock On the optimal control of {SIR} model with {Erlang-distributed} infectious period: isolation strategies.
\newblock {\em Journal of Mathematical Biology}, 83(4):36, Sep 2021.
\newblock \href {https://doi.org/10.1007/s00285-021-01668-1} {\path{doi:10.1007/s00285-021-01668-1}}.

\bibitem{BUONOMO2014}
B.~Buonomo, D.~Lacitignola, and C.~Vargas-De-León.
\newblock Qualitative analysis and optimal control of an epidemic model with vaccination and treatment.
\newblock {\em Mathematics and Computers in Simulation}, 100:88--102, 2014.
\newblock \href {https://doi.org/10.1016/j.matcom.2013.11.005} {\path{doi:10.1016/j.matcom.2013.11.005}}.

\bibitem{MPEuler}
H.~Burchard, E.~Deleersnijder, and A.~Meister.
\newblock A high-order conservative {Patankar-type} discretisation for stiff systems of production–destruction equations.
\newblock {\em Applied Numerical Mathematics}, 47(1):1--30, 2003.
\newblock \href {https://doi.org/10.1016/S0168-9274(03)00101-6} {\path{doi:10.1016/S0168-9274(03)00101-6}}.

\bibitem{Cacace_Oliviero}
S.~Cacace and A.~Oliviero.
\newblock Reliable optimal controls for seir models in epidemiology.
\newblock {\em Mathematics and Computers in Simulation}, 223:523--542, 2024.
\newblock \href {https://doi.org/10.1016/j.matcom.2024.04.034} {\path{doi:10.1016/j.matcom.2024.04.034}}.

\bibitem{Epidemic_PDS_1}
E.~Lima Campos, R.~Penha Cysne, A.~L. Madureira, and G.~L.Q. Mendes.
\newblock {Multi-generational SIR modeling: Determination of parameters, epidemiological forecasting and age-dependent vaccination policies}.
\newblock {\em Infectious Disease Modelling}, 6:751--765, 2021.
\newblock \href {https://doi.org/10.1016/j.idm.2021.05.003} {\path{doi:10.1016/j.idm.2021.05.003}}.

\bibitem{CAPASSO197843}
V.~Capasso and G.~Serio.
\newblock {A generalization of the Kermack-McKendrick deterministic epidemic model}.
\newblock {\em Mathematical Biosciences}, 42(1):43--61, 1978.
\newblock \href {https://doi.org/10.1016/0025-5564(78)90006-8} {\path{doi:10.1016/0025-5564(78)90006-8}}.

\bibitem{Epidemic_PDS_2}
A.~Cardone, P.~D. de~Alba, and B.~Paternoster.
\newblock {Analytical Properties and Numerical Preservation of an Age-Group Susceptible-Infected-Recovered Model: Application to the Diffusion of Information}.
\newblock {\em Journal of Computational and Nonlinear Dynamics}, 19(6):061006, 05 2024.
\newblock \href {https://doi.org/10.1115/1.4065437} {\path{doi:10.1115/1.4065437}}.

\bibitem{CFPS24}
E.~Carlini, R.~Ferretti, S.~Preda, and M.~Semplice.
\newblock A {CWENO} large time-step scheme for {Hamilton}--{Jacobi} equations.
\newblock {\em arXiv preprint arXiv:2402.15367}, 2024.

\bibitem{Cesari1965existence}
L.~Cesari.
\newblock An existence theorem in problems of optimal control.
\newblock {\em Journal of the Society for Industrial and Applied Mathematics, Series A: Control}, 3(1):7--22, 1965.

\bibitem{CdS}
M.~Ceseri, R.~Natalini, and M.~Pezzella.
\newblock {An Integro-differential Model of Cadmium Yellow Photodegradation}.
\newblock {\em arXiv}, 2024.
\newblock \href {https://doi.org/10.48550/arXiv.2411.06997} {\path{doi:10.48550/arXiv.2411.06997}}.

\bibitem{Reinfection}
Y.~Chen, W.~Zhu, X.~Han, M.~Chen, X.~Li, H.~Huang, M.~Zhang, R.~Wei, H.~Zhang, C.~Yang, and T.~Zhang.
\newblock {How does the SARS-CoV-2 reinfection rate change over time? The global evidence from systematic review and meta-analysis}.
\newblock {\em BMC Infectious Diseases}, 24(1):339, Mar 2024.
\newblock \href {https://doi.org/10.1186/s12879-024-09225-z} {\path{doi:10.1186/s12879-024-09225-z}}.

\bibitem{coppel1965stability}
W.A. Coppel.
\newblock {\em Stability and Asymptotic Behaviour of Differential Equations}.
\newblock Heath and mathematical monographs. HEATH AND COMP., 1965.

\bibitem{CrandallEvansLions1984}
M.~G. Crandall, L.~C. Evans, and P.-L. Lions.
\newblock Some properties of viscosity solutions of {Hamilton}--{Jacobi} equations.
\newblock {\em Transactions of the American Mathematical Society}, 282(2):487--502, 1984.
\newblock \href {https://doi.org/10.1090/S0002-9947-1984-0732102-X} {\path{doi:10.1090/S0002-9947-1984-0732102-X}}.

\bibitem{CrandallLions1983}
M.~G. Crandall and P.-L. Lions.
\newblock Viscosity solutions of {Hamilton}--{Jacobi} equations.
\newblock {\em Transactions of the American mathematical society}, 277(1):1--42, 1983.
\newblock \href {https://doi.org/10.1090/S0002-9947-1983-0690039-8} {\path{doi:10.1090/S0002-9947-1983-0690039-8}}.

\bibitem{CM10}
E.~Cristiani and P.~Martinon.
\newblock Initialization of the shooting method via the {Hamilton}--{Jacobi}--{Bellman} approach.
\newblock {\em Journal of Optimization Theory and Applications}, 146(2):321--346, 2010.
\newblock \href {https://doi.org/10.1007/s10957-010-9649-6} {\path{doi:10.1007/s10957-010-9649-6}}.

\bibitem{Dimitrov2}
D.~T.~Dimitrov D.~T.~Wood and H.~V. Kojouharov.
\newblock A nonstandard finite difference method for n-dimensional productive-destructive systems.
\newblock {\em Journal of Difference Equations and Applications}, 21(3):240--254, 2015.
\newblock \href {https://doi.org/10.1080/10236198.2014.997228} {\path{doi:10.1080/10236198.2014.997228}}.

\bibitem{DiekmannHeesterbeek2000}
O.~Diekmann and J.~A.~P. Heesterbeek.
\newblock {\em Mathematical epidemiology of infectious diseases: model building, analysis and interpretation}.
\newblock Wiley series in mathematical and computational biology. John Wiley \& Sons, United States, 2000.

\bibitem{Dimitrov1}
D.~T. Dimitrov and H.~V. Kojouharov.
\newblock Dynamically consistent numerical methods for general productive-destructive systems.
\newblock {\em Journal of Difference Equations and Applications}, 17(12):1721--1736, 2011.
\newblock \href {https://doi.org/10.1080/10236191003781947} {\path{doi:10.1080/10236191003781947}}.

\bibitem{dontchev2020approximating}
A.~L. Dontchev, I.~V. Kolmanovsky, M.~I. Krastanov, V.~M. Veliov, and P.~T. Vuong.
\newblock Approximating optimal finite horizon feedback by model predictive control.
\newblock {\em Systems \& Control Letters}, 139:104666, 2020.

\bibitem{DOnofrio}
A.~d’Onofrio.
\newblock {Vaccination policies and nonlinear force of infection: generalization of an observation by Alexander and Moghadas (2004)}.
\newblock {\em Applied Mathematics and Computation}, 168(1):613--622, 2005.
\newblock \href {https://doi.org/10.1016/j.amc.2004.09.013} {\path{doi:10.1016/j.amc.2004.09.013}}.

\bibitem{Evans1998}
L.C. Evans.
\newblock {\em Partial Differential Equations}.
\newblock Graduate studies in mathematics. American Mathematical Society, 1998.

\bibitem{Falcone_Ferretti_High}
M.~Falcone and R.~Ferretti.
\newblock {Discrete time high-order schemes for viscosity solutions of Hamilton--Jacobi--Bellman equations}.
\newblock {\em Numerische Mathematik}, 67(3):315--344, Apr 1994.
\newblock \href {https://doi.org/10.1007/s002110050031} {\path{doi:10.1007/s002110050031}}.

\bibitem{FalconeFerretti}
M.~Falcone and R.~Ferretti.
\newblock {\em Semi-Lagrangian Approximation Schemes for Linear and {Hamilton}--{Jacobi} Equations}.
\newblock Society for Industrial and Applied Mathematics, 2013.
\newblock \href {https://doi.org/10.1137/1.9781611973051} {\path{doi:10.1137/1.9781611973051}}.

\bibitem{FalconeFerretti2016}
M.~Falcone and R.~Ferretti.
\newblock Numerical methods for {Hamilton}--{Jacobi} type equations.
\newblock In R{\'e}mi Abgrall and Chi-Wang Shu, editors, {\em Handbook of Numerical Methods for Hyperbolic Problems}, volume~17 of {\em Handbook of Numerical Analysis}, pages 603--626. Elsevier, 2016.

\bibitem{FalconeGiorgi1999}
M.~Falcone and T.~Giorgi.
\newblock An approximation scheme for evolutive {Hamilton--Jacobi} equations.
\newblock In William~M. McEneaney, G.~George Yin, and Qing Zhang, editors, {\em Stochastic Analysis, Control, Optimization and Applications: A Volume in Honor of W.H. Fleming}, pages 289--303. Birkh{\"a}user Boston, Boston, MA, 1999.

\bibitem{FERNANDEZVILLAVERDE2022104318}
J.~Fernández-Villaverde and C.~I. Jones.
\newblock {Estimating and simulating a SIRD Model of COVID-19 for many countries, states, and cities}.
\newblock {\em Journal of Economic Dynamics and Control}, 140:104318, 2022.
\newblock Covid-19 Economics.
\newblock \href {https://doi.org/10.1016/j.jedc.2022.104318} {\path{doi:10.1016/j.jedc.2022.104318}}.

\bibitem{Formaggia_Scotti}
L.~Formaggia and A.~Scotti.
\newblock {Positivity and Conservation Properties of Some Integration Schemes for Mass Action Kinetics}.
\newblock {\em SIAM Journal on Numerical Analysis}, 49(3):1267--1288, 2011.
\newblock \href {https://doi.org/10.1137/100789592} {\path{doi:10.1137/100789592}}.

\bibitem{Bacteria_PDS}
I.~Hense and A.~Beckmann.
\newblock The representation of cyanobacteria life cycle processes in aquatic ecosystem models.
\newblock {\em Ecological Modelling}, 221(19):2330--2338, 2010.
\newblock \href {https://doi.org/10.1016/j.ecolmodel.2010.06.014} {\path{doi:10.1016/j.ecolmodel.2010.06.014}}.

\bibitem{Higham}
D.~J. Higham.
\newblock Modeling and simulating chemical reactions.
\newblock {\em SIAM Review}, 50(2):347--368, 2008.
\newblock \href {https://doi.org/10.1137/060666457} {\path{doi:10.1137/060666457}}.

\bibitem{Horn_2012}
R.~A. Horn and C.~R. Johnson.
\newblock {\em Matrix Analysis}.
\newblock Cambridge University Press, 2012.
\newblock \href {https://doi.org/10.1017/cbo9781139020411} {\path{doi:10.1017/cbo9781139020411}}.

\bibitem{Huang2019}
J.~Huang and C.-W. Shu.
\newblock {Positivity-Preserving Time Discretizations for Production--Destruction Equations with Applications to Non-equilibrium Flows}.
\newblock {\em Journal of Scientific Computing}, 78(3):1811--1839, Mar 2019.
\newblock \href {https://doi.org/10.1007/s10915-018-0852-1} {\path{doi:10.1007/s10915-018-0852-1}}.

\bibitem{Huang2019II}
J.~Huang, W.~Zhao, and C.-W. Shu.
\newblock {A Third-Order Unconditionally Positivity-Preserving Scheme for Production--Destruction Equations with Applications to Non-equilibrium Flows}.
\newblock {\em Journal of Scientific Computing}, 79(2):1015--1056, May 2019.
\newblock \href {https://doi.org/10.1007/s10915-018-0881-9} {\path{doi:10.1007/s10915-018-0881-9}}.

\bibitem{IzginSSPMPRK}
{Huang, J.}, {Izgin, T.}, {Kopecz, S.}, {Meister, A.}, and {Shu, C.-W.}
\newblock {On the stability of strong-stability-preserving modified Patankar--Runge--Kutta schemes}.
\newblock {\em ESAIM: M2AN}, 57(2):1063--1086, 2023.
\newblock \href {https://doi.org/10.1051/m2an/2023005} {\path{doi:10.1051/m2an/2023005}}.

\bibitem{Ingalls2013}
B.~P. Ingalls.
\newblock {\em Mathematical modeling in systems biology}.
\newblock MIT Press, London, England, 2013.

\bibitem{GECO2}
T.~Izgin, S.~Kopecz, A.~Martiradonna, and A.~Meister.
\newblock On the dynamics of first and second order {GeCo} and {gBBKS} schemes.
\newblock {\em Applied Numerical Mathematics}, 193:43--66, 2023.
\newblock \href {https://doi.org/10.1016/j.apnum.2023.07.014} {\path{doi:10.1016/j.apnum.2023.07.014}}.

\bibitem{Izgin1}
T.~Izgin, S.~Kopecz, and A.~Meister.
\newblock {Recent Developments in the Field of Modified Patankar--Runge--Kutta methods}.
\newblock {\em PAMM}, 21(1):e202100027, 2021.
\newblock \href {https://doi.org/10.1002/pamm.202100027} {\path{doi:10.1002/pamm.202100027}}.

\bibitem{Izgin3}
T.~Izgin, S.~Kopecz, and A.~Meister.
\newblock {On the Stability of Unconditionally Positive and Linear Invariants Preserving Time Integration Schemes}.
\newblock {\em SIAM Journal on Numerical Analysis}, 60(6):3029--3051, 2022.
\newblock \href {https://doi.org/10.1137/22M1480318} {\path{doi:10.1137/22M1480318}}.

\bibitem{Izgin2}
{Izgin, T.}, {Kopecz, S.}, and {Meister, A.}
\newblock {On Lyapunov stability of positive and conservative time integrators and application to second order modified Patankar--Runge--Kutta schemes}.
\newblock {\em ESAIM: M2AN}, 56(3):1053--1080, 2022.
\newblock \href {https://doi.org/10.1051/m2an/2022031} {\path{doi:10.1051/m2an/2022031}}.

\bibitem{MPLM}
G.~Izzo, E.~Messina, M.~Pezzella, and Antonia Vecchio.
\newblock {Modified Patankar Linear Multistep methods for production--destruction systems}, 2024.
\newblock \href {https://arxiv.org/abs/2407.12540} {\path{arXiv:2407.12540}}, \href {https://doi.org/10.48550/arXiv.2407.12540} {\path{doi:10.48550/arXiv.2407.12540}}.

\bibitem{Kopecz2018}
S.~Kopecz and A.~Meister.
\newblock On order conditions for modified {P}atankar–{R}unge–{K}utta schemes.
\newblock {\em Appl. Numer. Math.}, 123:159--179, 2018.
\newblock \href {https://doi.org/10.1016/j.apnum.2017.09.004} {\path{doi:10.1016/j.apnum.2017.09.004}}.

\bibitem{Kopecz2018second}
S.~Kopecz and A.~Meister.
\newblock Unconditionally positive and conservative third order modified {Patankar--Runge--Kutta} discretizations of production--destruction systems.
\newblock {\em BIT Numerical Mathematics}, 58(3):691--728, Sep 2018.
\newblock \href {https://doi.org/10.1007/s10543-018-0705-1} {\path{doi:10.1007/s10543-018-0705-1}}.

\bibitem{Kopecz2019}
S.~Kopecz and A.~Meister.
\newblock {On the existence of three-stage third-order modified Patankar--Runge--Kutta schemes}.
\newblock {\em Numerical Algorithms}, 81(4):1473--1484, Aug 2019.
\newblock \href {https://doi.org/10.1007/s11075-019-00680-3} {\path{doi:10.1007/s11075-019-00680-3}}.

\bibitem{Michaelis2}
S.~C. Kou, Binny~J. Cherayil, Wei Min, Brian~P. English, and X.~Sunney Xie.
\newblock {Single-Molecule Michaelis--Menten Equations}.
\newblock {\em The Journal of Physical Chemistry B}, 109(41):19068--19081, 2005.
\newblock PMID: 16853459.
\newblock \href {https://doi.org/10.1021/jp051490q} {\path{doi:10.1021/jp051490q}}.

\bibitem{Laidler}
K.~J. Laidler.
\newblock A glossary of terms used in chemical kinetics, including reaction dynamics (iupac recommendations 1996).
\newblock {\em Pure and Applied Chemistry}, 68(1):149--192, 1996.
\newblock \href {https://doi.org/doi:10.1351/pac199668010149} {\path{doi:doi:10.1351/pac199668010149}}.

\bibitem{Liu1986}
W.-m. Liu, S.~A. Levin, and Y.~Iwasa.
\newblock Influence of nonlinear incidence rates upon the behavior of {SIRS} epidemiological models.
\newblock {\em Journal of Mathematical Biology}, 23(2):187--204, 1986.
\newblock \href {https://doi.org/10.1007/BF00276956} {\path{doi:10.1007/BF00276956}}.

\bibitem{GECO1}
A.~Martiradonna, G.~Colonna, and F.~Diele.
\newblock {GeCo: Geometric Conservative nonstandard schemes for biochemical systems}.
\newblock {\em Applied Numerical Mathematics}, 155:38--57, 2020.
\newblock Structural Dynamical Systems: Computational Aspects held in Monopoli (Italy) on June 12-15, 2018.
\newblock \href {https://doi.org/10.1016/j.apnum.2019.12.004} {\path{doi:10.1016/j.apnum.2019.12.004}}.

\bibitem{MMH12}
M.~McAsey, L.~Mou, and W.~Han.
\newblock Convergence of the forward-backward sweep method in optimal control.
\newblock {\em Computational Optimization and Applications}, 53(1):207--226, 2012.
\newblock \href {https://doi.org/10.1007/s10589-011-9454-7} {\path{doi:10.1007/s10589-011-9454-7}}.

\bibitem{MPVNSFD}
E.~Messina, M.~Pezzella, and A.~Vecchio.
\newblock A non-standard numerical scheme for an age-of-infection epidemic model.
\newblock {\em Journal of Computational Dynamics}, 9(2):239--252, 2022.
\newblock \href {https://doi.org/10.3934/jcd.2021029} {\path{doi:10.3934/jcd.2021029}}.

\bibitem{MPV_Axioms}
E.~Messina, M.~Pezzella, and A.~Vecchio.
\newblock {Positive Numerical Approximation of Integro-Differential Epidemic Model}.
\newblock {\em Axioms}, 11(2), 2022.
\newblock \href {https://doi.org/10.3390/axioms11020069} {\path{doi:10.3390/axioms11020069}}.

\bibitem{MPV_Mixing}
E.~Messina, M.~Pezzella, and A.~Vecchio.
\newblock A long-time behavior preserving numerical scheme for age-of-infection epidemic models with heterogeneous mixing.
\newblock {\em Applied Numerical Mathematics}, 200:344--357, 2024.
\newblock New Trends in Approximation Methods and Numerical Analysis {(FAATNA20$>$22)}.
\newblock \href {https://doi.org/10.1016/j.apnum.2023.04.009} {\path{doi:10.1016/j.apnum.2023.04.009}}.

\bibitem{Mickens_root}
R.~E. Mickens.
\newblock {A SIR-model with square-root dynamics: An NSFD scheme}.
\newblock {\em Journal of Difference Equations and Applications}, 16(2-3):209--216, 2010.
\newblock \href {https://doi.org/10.1080/10236190802495311} {\path{doi:10.1080/10236190802495311}}.

\bibitem{CudaGuide}
NVIDIA.
\newblock Cuda {C++} programming guide.
\newblock 2024.
\newblock URL: \url{https://docs.nvidia.com/cuda/pdf/CUDA_C_Programming_Guide.pdf}.

\bibitem{Torlo2020}
P.~{\"O}ffner and D.~Torlo.
\newblock Arbitrary high-order, conservative and positivity preserving {Patankar-type} deferred correction schemes.
\newblock {\em Applied Numerical Mathematics}, 153:15--34, 2020.
\newblock \href {https://doi.org/10.1016/j.apnum.2020.01.025} {\path{doi:10.1016/j.apnum.2020.01.025}}.

\bibitem{patankar1980numerical}
S.~V. Patankar.
\newblock {\em Numerical heat transfer and fluid flow}.
\newblock Series on Computational Methods in Mechanics and Thermal Science. Hemisphere Publishing Corporation (CRC Press, Taylor \& Francis Group), 1980.

\bibitem{Mario_Pezzella}
M.~Pezzella.
\newblock High order positivity-preserving numerical methods for a non-local photochemical model.
\newblock {\em arXiv}, 2025.
\newblock \href {https://doi.org/10.48550/arXiv.2501.04573} {\path{doi:10.48550/arXiv.2501.04573}}.

\bibitem{roxin1962existence}
E.~Roxin.
\newblock The existence of optimal controls.
\newblock {\em Michigan Mathematical Journal}, 9(2):109--119, 1962.

\bibitem{SR18}
R.~A. Sheldon and J.~M. Woodley.
\newblock {Role of Biocatalysis in Sustainable Chemistry}.
\newblock {\em Chemical Reviews}, 118(2):801--838, 2018.
\newblock \href {https://doi.org/10.1021/acs.chemrev.7b00203} {\path{doi:10.1021/acs.chemrev.7b00203}}.

\bibitem{SSGK23}
M.~Sperl, L.~Saluzzi, L.~Gr{\"u}ne, and D.~Kalise.
\newblock Separable approximations of optimal value functions under a decaying sensitivity assumption.
\newblock In {\em 2023 62nd IEEE Conference on Decision and Control (CDC)}, pages 259--264. IEEE, 2023.

\bibitem{TORLO2022}
D.~Torlo, P.~{\"O}ffner, and H.~Ranocha.
\newblock Issues with positivity-preserving {P}atankar-type schemes.
\newblock {\em Appl. Numer. Math.}, 182:117--147, 2022.
\newblock \href {https://doi.org/10.1016/j.apnum.2022.07.014} {\path{doi:10.1016/j.apnum.2022.07.014}}.

\bibitem{vandenDriessche2000}
P.~van~den Driessche and J.~Watmough.
\newblock {A simple SIS epidemic model with a backward bifurcation}.
\newblock {\em Journal of Mathematical Biology}, 40(6):525--540, Jun 2000.
\newblock \href {https://doi.org/10.1007/s002850000032} {\path{doi:10.1007/s002850000032}}.

\bibitem{WALSH20101}
R.~Walsh, E.~Martin, and S.~Darvesh.
\newblock A method to describe enzyme-catalyzed reactions by combining steady state and time course enzyme kinetic parameters.
\newblock {\em Biochimica et Biophysica Acta (BBA) - General Subjects}, 1800(1):1--5, 2010.
\newblock \href {https://doi.org/10.1016/j.bbagen.2009.10.007} {\path{doi:10.1016/j.bbagen.2009.10.007}}.

\bibitem{Zafar}
Z.~A. Zafar, K.~Rehan, M.~Mushtaq, and M.~Rafiq.
\newblock {Numerical Modeling for Nonlinear Biochemical Reaction Networks}.
\newblock {\em Iranian Journal of Mathematical Chemistry}, 8(4):413--423, 2017.
\newblock \href {https://doi.org/10.22052/ijmc.2017.47506.1170} {\path{doi:10.22052/ijmc.2017.47506.1170}}.

\bibitem{Zhang2024-br}
L.~Zhang, J.~Peng, Y.~Ge, H.~Li, and Y.~Tang.
\newblock {High-Accuracy} {Positivity-Preserving} {Finite Difference Approximations of the Chemotaxis Model for Tumor Invasion}.
\newblock {\em Journal of Computational Biology}, October 2024.
\newblock \href {https://doi.org/10.1089/cmb.2023.0316} {\path{doi:10.1089/cmb.2023.0316}}.

\bibitem{Zhu2024}
F.~Zhu, J.~Huang, and Y.~Yang.
\newblock {Bound-Preserving Discontinuous Galerkin Methods with Modified Patankar Time Integrations for Chemical Reacting Flows}.
\newblock {\em Communications on Applied Mathematics and Computation}, 6(1):190--217, Mar 2024.
\newblock \href {https://doi.org/10.1007/s42967-022-00231-z} {\path{doi:10.1007/s42967-022-00231-z}}.

\end{thebibliography}

%
%

\end{document}